\documentclass[review]{siamart251216}
\usepackage{amssymb,amscd,amsxtra,calc}
\usepackage{cmmib57}
\usepackage{multirow}
\usepackage{mathrsfs}
\usepackage{graphicx}
\usepackage{bm}
\usepackage{amsmath}
\usepackage{booktabs}
\usepackage{amssymb}
\usepackage{graphicx}
\usepackage{accents}
\usepackage{tikz} 
\usepackage{tikz-cd}
\usepackage{pgfplots}
\pgfplotsset{compat=1.18}
\usepgfplotslibrary{groupplots}
\definecolor{seaborngreen}{rgb}{0.3333333333333333, 0.6588235294117647, 0.40784313725490196}
\definecolor{seabornblue}{rgb}{0.2980392156862745, 0.4470588235294118, 0.6901960784313725}
\definecolor{seabornred}{rgb}{0.7686274509803922, 0.3058823529411765, 0.3215686274509804}
\usepackage{stmaryrd}
\usepackage{comment}
\usepackage{colonequals}
\usepackage{subcaption}
\usepackage{algorithmic}
\usepackage{algorithm}
\usepackage{graphicx}
\allowdisplaybreaks[3]
\renewcommand{\arraystretch}{1.3}

\hypersetup{
    colorlinks,
    allcolors={green!50!black},
    urlcolor={red!50!black}
}

\usepackage{colonequals}

\newsiamthm{example}{Example}
\newsiamremark{remark}{Remark}
\crefname{example}{Example}{Examples}
\Crefname{example}{Example}{Examples}
\crefname{remark}{Remark}{Remarks}
\Crefname{remark}{Remark}{Remarks}
\newcommand{\elemma}{\end{lemma}}
\newcommand{\blemma}{\begin{lemma}}
\theoremstyle{remark}

\numberwithin{equation}{section}

\newcommand{\btheorem}{\begin{theorem}}
	\newcommand{\etheorem}{\end{theorem}}
\newcommand{\bproposition}{\begin{prop}}
	\newcommand{\eproposition}{\end{prop}}
\newcommand{\bdefinition}{\begin{defn}}
	\newcommand{\edefinition}{\end{defn}}
\newcommand{\bcorollary}{\begin{cor}}
	\newcommand{\ecorollary}{\end{cor}}
\newcommand{\bproof}{\begin{proof}}
	\newcommand{\eproof}{\end{proof}}
\newcommand{\bremark}{\begin{remark}}	
	\newcommand{\eremark}{\end{remark}}
\newcommand{\eexample}{\end{example}}
\newcommand{\bexample}{\begin{example}}

\newcommand{\bld}{\boldsymbol}
\newcommand{\lb}{\left(}
\newcommand{\rb}{\right)}
\newcommand{\lmb}{\left[}
\newcommand{\rmb}{\right]}

\newcommand{\upd}{\mathrm{d }}
\newcommand{\p}{\partial}

\renewcommand{\S}{\mathbb{S}}

\newcommand{\bomega}{\boldsymbol{\omega}}

\newcommand{\bX}{\boldsymbol{X}}
\newcommand{\bY}{\boldsymbol{Y}}
\newcommand{\bU}{\boldsymbol{U}}

\newcommand{\bphi}{\boldsymbol{\varphi}}
\newcommand{\bpsi}{\boldsymbol{\psi}}

\newcommand{\mH}{\mathcal{H}}

\renewcommand\bX{{\boldsymbol{X}}}

\renewcommand{\S}{\mathbb{S}}

\nolinenumbers

\headers{IMEX Runge-Kutta schemes for geometric flows}%
        {W.~Jiang, C.~Su, and K.~Tang}
\title{Implict-Explicit Runge-Kutta schemes for geometric flows\thanks{Submitted to the editors DATE.
\funding{This work is supported by the NSFC 12671530 and 12271414 (Wei Jiang), and NSFC 12201342 (Chunmei Su).}}}
\author{Wei Jiang\thanks{School of Mathematics and Statistics, Hubei Key Laboratory of Computational Science, Wuhan University, Wuhan 430072, China (\email{jiangwei1007@whu.edu.cn}).}
\and Chunmei Su\thanks{Yau Mathematical Sciences Center, Tsinghua University, Beijing, 100084, China
		(\email{sucm@tsinghua.edu.cn}).} \and Kaike Thang\thanks{Department of Mathematics, Tsinghua University, Beijing, 100084, China
  (\email{tkk20@mails.tsinghua.edu.cn}).}}

\begin{document}

\maketitle

\begin{abstract}
We propose a one-step implicit-explicit (IMEX) Runge-Kutta framework for high-order time integration of geometric flows. Starting from the Dziuk and Deckelnick-Dziuk parametric finite element formulations of curve shortening flow, we recast the spatially semidiscrete equations as systems of ordinary differential equations and construct second- and third-order schemes. Each stage requires only a linear solve. The key ingredient is a reduced-coefficient structure that allows the variational problem at each stage to be posed on the geometry generated by the preceding stage. The Runge-Kutta stages thus supply the required prediction geometry without a separate prediction procedure. The resulting schemes require no additional starting values. The framework also extends directly to the Barrett-Garcke-N\"urnberg formulations of mean curvature flow and surface diffusion for curves and surfaces. Numerical experiments support the designed orders of temporal convergence and indicate that suitable coefficient choices yield mesh quality comparable to that of the corresponding first-order methods.
\end{abstract}

\begin{keywords}
 geometric flow, parametric finite element method, Runge-Kutta method, mesh quality, Dziuk scheme, BGN formulation.
\end{keywords}

\begin{MSCcodes}
  65M60, 65M12, 53E10, 35K93
\end{MSCcodes}

\section{Introduction}\label{sec:intro}

Geometric flows are evolution equations in which the motion of geometric objects is governed by intrinsic or extrinsic geometric quantities. They arise in a broad range of scientific and engineering applications, including grain-boundary motion in materials science \cite{mullins1956two}, solid-state dewetting of thin films \cite{bao2017parametric,zhao2020parametric}, image processing and computer vision \cite{cao2003geometric}, and biomembrane modeling \cite{du2005phase}. Their numerical approximation is also closely connected to threshold-dynamics algorithms \cite{esedoglu2008threshold,ruuth1999convolution}. Accurate and robust numerical methods are therefore essential for capturing the evolving geometry and gaining insight into the qualitative behavior of these flows.

The parametric finite element method (PFEM) is a widely used approach to the numerical approximation of geometric flows. Originally introduced by Dziuk \cite{dziuk1990algorithm} for mean curvature flow of surfaces in three-dimensional space, it has since been extended to a broad class of geometric evolution problems. Representative examples include mean curvature flow \cite{dziuk1994convergence,kovacs2019convergent,kovacs2020convergent}, surface diffusion flow \cite{bao2021structure,jiang2021perimeter}, Willmore flow \cite{dziuk2002evolution,dziuk2008computational,kovacs2021convergent,bao2025energy}, and anisotropic geometric flows \cite{dziuk1999discrete,barrett2008numerical,bao2023symmetrized,deckelnick2025error}. For comprehensive accounts of parametric finite element methods and their applications to geometric evolution problems, we refer the reader to \cite{deckelnick2005computation,barrett2020parametric,bansch2023interfaces}.

Maintaining mesh quality throughout the evolution remains a persistent challenge in PFEMs. In particular, node clustering and excessive element distortion can substantially degrade numerical accuracy and may ultimately lead to computational instability. Various strategies have been proposed to address this issue, including artificial mesh regularization \cite{bansch2005finite}, prescribed tangential velocities \cite{sevcovic2001evolution,mikula2004direct}, DeTurck-type harmonic-map techniques \cite{Deckelnick1995,m2017approximations}, and discrete harmonic-map reparametrization \cite{steinhilber2014numerical}. Hu and Li \cite{hu2022evolving} developed an evolving surface finite element method that incorporates an artificial tangential velocity to improve mesh quality in curvature-driven flows. More recently, Duan and Li \cite{duan2024new} proposed a class of PFEMs in which the tangential velocity is determined by a harmonic map between a fixed reference surface and the evolving surface.

An effective alternative is the parametric finite element approach developed by Barrett, Garcke, and N{\"u}rnberg (BGN) \cite{barrett2007parametric}. The BGN formulation approximates the prescribed normal evolution while allowing tangential motion to arise naturally from the discretization. This tangential motion can substantially improve the mesh distribution without the need to prescribe a tangential velocity explicitly. BGN-type schemes have been applied to a wide range of geometric flows; see, for example, \cite{barrett2007variational,barrett2010numerical,bao2021structure,bao2023symmetrized}. Recently, Gao et al. \cite{gao2026energy} investigated the connection between the tangential motion induced by the BGN formulation and minimal deformation. Building on this connection, they developed the BGN-MDR method, which is a linear, energy-stable scheme with a nontrivial tangential velocity.

Many classical PFEMs employ first-order time discretizations. Achieving higher-order temporal accuracy is more challenging because the geometry, the finite element space, and the mesh distribution evolve simultaneously. Kov{\'a}cs, Li, and Lubich established error estimates for high-order evolving finite element methods for several geometric flows by reformulating the semidiscrete equations as differential-algebraic systems and applying backward differentiation formula (BDF) time discretizations \cite{kovacs2019convergent,kovacs2020convergent,kovacs2021convergent}. Duan and Li \cite{duan2024new} also considered a second-order time discretization.

For BGN-type schemes, Jiang, Su, and Zhang \cite{jiang2024second} constructed a second-order method based on a semi-implicit Crank--Nicolson leap-frog discretization. Higher-order BDF schemes were subsequently developed in \cite{jiang2024stable}. A key ingredient in these BGN-based methods is the construction of suitable predictor curves, which play an important role in ensuring temporal consistency and controlling mesh redistribution. Predictor--corrector strategies offer another route to second-order temporal accuracy \cite{jiang2025predictor}, and related ideas have recently been applied to curve shortening flow and curve diffusion \cite{deckelnick2026second}. More recently, Zhang, Andrews, and Farrell \cite{zhang2026arbitrary} developed structure-preserving discretizations of arbitrarily high order in both space and time for curvature flows. Their approach combines suitable auxiliary variables with a continuous Petrov--Galerkin discretization in time and incorporates the minimal deformation rate strategy to maintain good mesh quality.

Many of the high-order PFEMs discussed above employ multistep time discretizations. Alternatives include the predictor--corrector method \cite{jiang2025predictor} and the recent space--time Petrov--Galerkin approach \cite{zhang2026arbitrary}. A practical drawback of multistep methods is that they require additional starting values of sufficient accuracy, typically generated through a separate initialization procedure using auxiliary schemes. For BGN-type formulations, high-order time discretization presents an additional geometric challenge. Our numerical experience suggests that the choice of predictor curves used in the variational formulation is crucial to retaining the favorable mesh-distribution properties of the BGN formulation~\cite{jiang2024stable,jiang2025predictor}. An unsuitable choice may lead to severe mesh distortion and ultimately cause computational breakdown. High-order BGN-type schemes must therefore combine high-order temporal accuracy with effective tangential mesh redistribution. The need to satisfy both requirements makes the design of such schemes particularly challenging.

Runge--Kutta (RK) methods provide a natural alternative to multistep time discretizations. They do not require multiple starting values, and their multistage structure naturally supplies intermediate geometric information within each time step. This feature is particularly attractive for geometric flows, where the geometry used at intermediate stages can strongly affect mesh redistribution. Consider, for example, an ODE system of the form
\begin{align}
\bld{U}_t = \bld{f}(t,\bld{U})
+ G(t,\bld{U})\bld{U}.
\label{IMEX-RK: model problem}
\end{align}
For stiff problems, fully explicit RK methods generally impose restrictive time-step constraints. Implicit-explicit (IMEX) RK methods alleviate this restriction by treating different components of the system with different levels of implicitness. When the stiff operator is linear, this strategy leads to efficient stage equations and has been successfully applied to the Allen--Cahn and Cahn--Hilliard equations \cite{song2016energy,fu2022energy,fu2024energy}, as well as to hyperbolic systems \cite{pareschi2005implicit,boscarino2009class,boscarino2017unified}. The situation is more delicate when $G$ depends nonlinearly on $\bld{U}$. A fully implicit treatment of $G(t,\bld{U})\bld{U}$ would generally require the solution of a nonlinear system at every RK stage. Ding \cite{ding2025semi} recently proposed a class of RK methods that avoids such nonlinear solves by treating $G(t,\bld{U})\bld{U}$ semi-implicitly at the intermediate stages. This idea is particularly suitable for parametric geometric flows, where a carefully chosen stage geometry can simultaneously linearize the stage equations and provide the geometric information needed for mesh redistribution.

In this paper, we develop a high-order implicit--explicit Runge--Kutta (IMEX-RK) framework for the parametric finite element approximation of geometric flows. Starting from Dziuk's semidiscrete formulation of curve shortening flow (CSF), we recast the resulting finite-dimensional equations as an ordinary differential equation (ODE) system suitable for IMEX-RK discretization. This construction yields $s$-stage methods of order $s$ for $s=2,3$, with each stage requiring only a linear solve. A key feature is that the intermediate RK stages themselves generate the predicted geometries needed at subsequent stages, eliminating the need for separate extrapolation or auxiliary predictor curves. These stage geometries play a role analogous to that of the predictor curves used in high-order BDF methods \cite{jiang2024stable}, but arise intrinsically from the RK discretization itself. {\it The main contributions of this work are summarized as follows}:
\smallskip
\begin{itemize}

\item We systematically construct a family of high-order IMEX-RK schemes for the PFEM approximation of CSF. For the orders considered, these schemes achieve $s$th-order temporal accuracy with $s$ RK stages, each requiring only a linear solve.
We exploit the intrinsic multistage structure of RK methods to generate the required geometric predictions. Each completed stage supplies the geometry needed at subsequent stages, eliminating the need for the separate predictor curves often required in high-order multistep discretizations.

\item Although initially developed for Dziuk's method, the IMEX-RK framework extends naturally to the BGN formulation. This extension demonstrates its applicability to different parametrization strategies and its compatibility with the intrinsic tangential mesh redistribution of BGN-type schemes.

\item The proposed IMEX-RK framework extends beyond CSF to other geometric flows, including surface diffusion. It also applies to the evolution of two-dimensional surfaces embedded in three-dimensional space.

\item The proposed high-order schemes retain the mesh-redistribution mechanism of the underlying first-order PFEM. Numerical experiments indicate that, with suitable RK coefficients, these schemes achieve high-order temporal accuracy while maintaining mesh quality comparable to that of the original method.

\item Numerical experiments on representative curve evolution problems support the expected orders of temporal convergence and demonstrate the accuracy, computational efficiency, and favorable mesh quality of the proposed methods.
\end{itemize}

\smallskip

The remainder of this paper is organized as follows. Section~2 develops the IMEX-RK discretization of Dziuk's semidiscrete scheme for CSF. Section~3 extends the construction to the BGN formulation. Section~4 presents numerical experiments that examine the accuracy, efficiency, and mesh behavior of the proposed high-order methods. Finally, Section~5 summarizes the main results and gives concluding remarks.

\section{Implicit-Explicit Runge-Kutta temporal discretization}

\smallskip
In this section, we develop the IMEX--RK temporal discretization for
parametric finite element approximations of geometric flows. We begin by
reviewing two classical PFEM formulations for CSF: Dziuk's formulation and the Deckelnick--Dziuk formulation. We then recast
their semidiscrete systems as finite-dimensional ODEs, which provides the
starting point for the high-order IMEX--RK schemes developed below.

\subsection{Dziuk-type PFEMs for CSF}
\smallskip
Let $\S^1=\mathbb{R}/(2\pi\mathbb{Z})$, represented by the interval $[0,2\pi]$ with its endpoints identified. Let $\bX:\S^1\to\mathbb{R}^2$ be a smooth immersion, and denote the associated closed curve by
\[
\Gamma[\bX]=\bX(\S^1).
\]
We assume that $\bX$ parametrizes $\Gamma[\bX]$ in the anticlockwise direction.
The arc-length derivative is defined by
\(
\p_s f=\frac{\p_\rho f}{|\p_\rho\bX|}.
\)
Let
\[
\bld{N}=-(\p_s\bX)^\perp,
\qquad
(a,b)^\perp=(-b,a),
\]
be the outward unit normal. With this convention, the curvature is given by
\[
\kappa=-\p_{ss}\bX\cdot\bld{N}.
\]

The length of $\Gamma$ is
\(
L(\Gamma)=\int_\Gamma 1\,\upd s. \)
Its $L^2$-gradient flow is the curve shortening flow (CSF)
\begin{equation}
\p_t\bX
=-\kappa\bld{N}
=
\p_{ss}\bX.
\label{IMEX-RK-CSF: CSF in geometric form}
\end{equation}
Since
\(
\p_{ss}\bX
=\frac{1}{|\p_\rho\bX|}
\p_\rho
\left(
\frac{\p_\rho\bX}{|\p_\rho\bX|}
\right),
\)
the evolution equation can be written in the reference coordinate $\rho$ as
\begin{equation}
\p_t\bX
=\frac{1}{|\p_\rho\bX|}
\p_\rho
\left(
\frac{\p_\rho\bX}{|\p_\rho\bX|}
\right).
\label{IMEX-RK-CSF: CSF in reference coordinates}
\end{equation}
Multiplying \eqref{IMEX-RK-CSF: CSF in reference coordinates} by
$|\p_\rho\bX|$, testing with
$\bphi\in H^1(\S^1,\mathbb{R}^2)$, and integrating by parts yield the weak
formulation
\begin{align}
\int_{\S^1}
|\p_\rho\bX|\,
\p_t\bX\cdot\bphi\,\upd\rho
+
\int_{\S^1}
\frac{\p_\rho\bX\cdot\p_\rho\bphi}
{|\p_\rho\bX|}
\,\upd\rho
=
0,
\qquad
\forall\,\bphi\in H^1(\S^1,\mathbb{R}^2).
\label{IMEX-RK-CSF: CSF weak form}
\end{align}

We next introduce the spatial finite element discretization. Let
\(
0=\rho_0<\rho_1<\cdots<\rho_N=2\pi
\) be a partition of $\S^1$, and set
\(I_j=[\rho_{j-1},\rho_j]\),
$h_j=\rho_j-\rho_{j-1}$,
$h=\max_{1\leq j\leq N}h_j$. The periodic piecewise linear finite element space is
\[V^h
:=\left\{
w\in C(\S^1,\mathbb{R})
:
w|_{I_j}\ \text{is linear},\ j=1,\ldots,N;\quad w(0)=w(2\pi)
\right\}.
\]

For
$\bX_h\in[V^h]^2$, define
\(\bX_j:=\bX_h(\rho_j)\), $\bld{h}_j:=\bX_j-\bX_{j-1}$,
and assume
\(
|\bld{h}_j|>0\),
$j=1,\ldots,N$. The polygonal curve is
\(
\Gamma_h=\bX_h(\S^1),
\) and its elementwise outward unit normal is
\(
\bld{N}_h|_{I_j}
=-\left(
\frac{\bld{h}_j}{|\bld{h}_j|}
\right)^\perp.
\)
For a piecewise smooth scalar or vector function $f$, we define
\[
\p_s f
=
\frac{\p_\rho f}{|\p_\rho\bX_h|}.
\]

For piecewise continuous scalar or vector functions $f$ and $g$, possibly
with jumps at the mesh nodes, the mass-lumped inner product on $\Gamma_h$
is defined by
\[
(f,g)^h_{\Gamma_h}
:=
\frac12
\sum_{j=1}^N
|\bld{h}_j|
\left[
(f\cdot g)(\rho_j^-)
+
(f\cdot g)(\rho_{j-1}^+)
\right],
\]
where
\(
u(\rho_j^\pm)
=
\lim\limits_{\varepsilon\rightarrow 0^+}
u(\rho_j\pm\varepsilon).
\)
We also use the standard $L^2$ inner product
\[
(f,g)_{\Gamma_h}
:=
\int_{\Gamma_h}f\cdot g\,\upd s.
\]
It can be easily seen that the mass-lumped inner product is an approximation of the $L^2$ inner product by using the trapezoid rule on each interval, which is exact if $f\cdot g$ is piecewise linear.

Dziuk's semidiscrete PFEM for CSF \cite{dziuk1994convergence} is obtained
by applying the above finite element discretization to
\eqref{IMEX-RK-CSF: CSF weak form}. Given an initial approximation
$\bX_h(\cdot,0)\in[V^h]^2$, find
$\bX_h(\cdot,t)\in[V^h]^2$ such that
\begin{align}
\int_{\S^1}
\operatorname{I}_h
\left(
\p_t\bX_h\cdot\bphi_h
\right)
|\p_\rho\bX_h|
\,\upd\rho
+
\int_{\S^1}
\frac{
\p_\rho\bX_h\cdot\p_\rho\bphi_h
}{
|\p_\rho\bX_h|
}
\,\upd\rho
=
0,
\qquad
\forall\,\bphi_h\in[V^h]^2,
\label{IMEX-RK-CSF: Dziuk's scheme}
\end{align}
where
\(
\operatorname{I}_h:C(\S^1,\mathbb{R})\rightarrow V^h
\)
denotes the nodal interpolation operator. Equivalently,
\eqref{IMEX-RK-CSF: Dziuk's scheme} can be written as
\begin{align}
\left(
\p_t\bX_h,\bphi_h
\right)^h_{\Gamma_h(t)}
+
\left(
\p_s\bX_h,\p_s\bphi_h
\right)_{\Gamma_h(t)}
=
0,
\qquad
\forall\,\bphi_h\in[V^h]^2,
\label{IMEX-RK-CSF: Dziuk compact}
\end{align}
with
\(
\Gamma_h(t)=\bX_h(\S^1,t).
\)

A semi-implicit Euler discretization of
\eqref{IMEX-RK-CSF: Dziuk compact} gives the classical first-order fully
discrete scheme. Given $\bX^m\in[V^h]^2$, find
$\bX^{m+1}\in[V^h]^2$ such that
\begin{align}
\left(
\frac{\bX^{m+1}-\bX^m}{\tau},
\bphi_h
\right)^h_{\Gamma^m}
+
\left(
\p_s\bX^{m+1},
\p_s\bphi_h
\right)_{\Gamma^m}
=
0,
\qquad
\forall\,\bphi_h\in[V^h]^2,
\label{IMEX-RK-CSF: Dziuk fully discrete}
\end{align}
where
\(
\Gamma^m=\bX^m(\S^1).
\)
Here and below, the geometric quantities appearing in the inner products
and arc-length derivatives are evaluated on the known curve $\Gamma^m$.
Although this scheme is simple and robust, its tangential mesh
redistribution is limited, and nonuniform meshes may develop during more
complicated evolutions.

Deckelnick and Dziuk  introduced a closely related
parametrization with stronger tangential motion \cite{Deckelnick1995}. They considered
\begin{align}
\p_t\bX
=
\frac{\p_{\rho\rho}\bX}{|\p_\rho\bX|^2}.
\label{IMEX-RK-CSF: Dziuk-Deckelnick formulation}
\end{align}
To compare this equation with CSF, observe that
\begin{align}
\frac{1}{|\p_\rho\bX|}
\p_\rho
\left(
\frac{\p_\rho\bX}{|\p_\rho\bX|}
\right)
=
\frac{\p_{\rho\rho}\bX}{|\p_\rho\bX|^2}
-
\frac{
\p_{\rho\rho}\bX\cdot\p_\rho\bX
}{
|\p_\rho\bX|^4
}
\p_\rho\bX.
\label{IMEX-RK-CSF: DD decomposition}
\end{align}
The second term on the right-hand side is tangential to the curve.
Consequently,
\eqref{IMEX-RK-CSF: Dziuk-Deckelnick formulation} and
\eqref{IMEX-RK-CSF: CSF in reference coordinates} have the same normal
velocity and therefore generate the same geometric evolution, but with
different parametrizations. The weak formulation of
\eqref{IMEX-RK-CSF: Dziuk-Deckelnick formulation} is
\begin{align}
\int_{\S^1}
|\p_\rho\bX|^2\,
\p_t\bX\cdot\bphi
\,\upd\rho
+
\int_{\S^1}
\p_\rho\bX\cdot\p_\rho\bphi
\,\upd\rho
=
0,
\qquad
\forall\,\bphi\in H^1(\S^1,\mathbb{R}^2).
\label{IMEX-RK-CSF: DD weak form}
\end{align}
Its mass-lumped semidiscrete approximation reads
\begin{align}
\int_{\S^1}
\operatorname{I}_h
\left(
\p_t\bX_h\cdot\bphi_h
\right)
|\p_\rho\bX_h|^2
\,\upd\rho
+
\int_{\S^1}
\p_\rho\bX_h\cdot\p_\rho\bphi_h
\,\upd\rho
=
0,
\qquad
\forall\,\bphi_h\in[V^h]^2.
\label{IMEX-RK-CSF: Dziuk-Deckelnick scheme}
\end{align}
or equivalently,
\begin{align*}
\left( |\p_\rho\bX_h|^2\p_t \bX_h,\bphi_h\right)^h_{\S^1} + \left( \p_\rho \bX_h,\p_\rho \bphi_h\right)_{\S^1} = 0,\quad \forall~ \bphi_h\in [V^h]^2,
\end{align*}
where $(f,g)_{\S^1}^h := \frac{1}{2}\sum\limits_{j=1}^N (\rho_j-\rho_{j-1})\lb (f\cdot g)(\rho_j^-) + (f\cdot g)(\rho_{j-1}^+)\rb
$, and $(\cdot,\cdot)_{\S^1}$ denotes the standard $L^2$ inner product on
$\S^1$.

A semi-implicit Euler discretization yields
\begin{align}
\left(
|\p_\rho\bX^m|^2
\frac{\bX^{m+1}-\bX^m}{\tau},
\bphi_h
\right)^h_{\S^1}
+
\left(
\p_\rho\bX^{m+1},
\p_\rho\bphi_h
\right)_{\S^1}
=
0,
\qquad
\forall\,\bphi_h\in[V^h]^2.
\label{IMEX-RK-CSF: Dziuk-Deckelnick fully discrete}
\end{align}
The schemes
\eqref{IMEX-RK-CSF: Dziuk fully discrete} and
\eqref{IMEX-RK-CSF: Dziuk-Deckelnick fully discrete}
are first-order accurate in time. Higher-order accuracy can be
obtained, for example, by replacing the backward Euler discretization with
a BDF formula. Such schemes are multistep methods and therefore require
sufficiently accurate starting values, which are generated by
additional auxiliary computations and might be very expensive to obtain sufficiently accurate starting values.

We seek single-step high-order temporal
discretizations of the semidiscrete systems
\eqref{IMEX-RK-CSF: Dziuk's scheme} and
\eqref{IMEX-RK-CSF: Dziuk-Deckelnick scheme}
based on Runge--Kutta methods. The multistage structure will allow the
intermediate RK solutions to provide the geometric information required
within each time step, without the separate start-up predictions required
by multistep discretizations.

\subsection{Construction of IMEX-RK methods}
We now construct high-order IMEX--RK temporal discretizations for the semidiscrete PFEMs
\eqref{IMEX-RK-CSF: Dziuk's scheme} and
\eqref{IMEX-RK-CSF: Dziuk-Deckelnick scheme}.
We first rewrite these semidiscrete formulations as finite-dimensional ODE systems.

For Dziuk's scheme, let
\[
\bU_h
=
(x_1,\ldots,x_N,y_1,\ldots,y_N)^T
\in\mathbb{R}^{2N},
\qquad
\bX_h(\rho_j)=(x_j,y_j)^T.
\]
Then \eqref{IMEX-RK-CSF: Dziuk's scheme} can be written as
\begin{align}
\frac{\upd\bU_h}{\upd t}
=
-A(\bU_h)^{-1}M(\bU_h)\bU_h,
\label{IMEX-RK-CSF: Dziuk and DD scheme in ODE form}
\end{align}
where
\[
A(\bU_h)
=
\operatorname{diag}\bigl(P(\bU_h),P(\bU_h)\bigr),
\qquad
M(\bU_h)
=
\operatorname{diag}\bigl(Q(\bU_h),Q(\bU_h)\bigr).
\]
The matrix $P(\bU_h)$ is diagonal, with
\(P_{jj}(\bU_h)
=\frac12
\left(
|\bld{h}_j|+|\bld{h}_{j+1}|
\right),
\)
and $Q(\bU_h)$ is the periodic stiffness tridiagonal matrix with entries
\[
Q_{jj}(\bU_h)
=
\frac{1}{|\bld{h}_j|}
+
\frac{1}{|\bld{h}_{j+1}|},
\qquad
Q_{j,j-1}(\bU_h)
=
-\frac{1}{|\bld{h}_j|},
\qquad
Q_{j,j+1}(\bU_h)
=
-\frac{1}{|\bld{h}_{j+1}|},
\]
where all indices are understood periodically modulo $N$.

For the Deckelnick--Dziuk formulation
\eqref{IMEX-RK-CSF: Dziuk-Deckelnick scheme}, we assume a uniform partition
with $h=2\pi/N$. It has the same ODE form
\eqref{IMEX-RK-CSF: Dziuk and DD scheme in ODE form}, but now
\[
P_{jj}(\bU_h)
=
\frac{1}{2h}
\left(
|\bld{h}_j|^2+|\bld{h}_{j+1}|^2
\right),\quad
Q_{jj}=\frac{2}{h},
\quad
Q_{j,j-1}=Q_{j,j+1}=-\frac{1}{h},
\]
again with periodic indexing.

Thus, both semidiscrete schemes can be written in the unified form
\begin{align}
\frac{\upd\bU}{\upd t}
=
\bld{F}(\bU)\bU,
\qquad
\bld{F}(\bU)
:=
-A(\bU)^{-1}M(\bU).
\label{IMEX-RK-CSF: unified ODE form}
\end{align}
In particular, \eqref{IMEX-RK-CSF: unified ODE form} is a special case of
\eqref{IMEX-RK: model problem} with
\(
\bld{f}(t,\bU)=0\),
\(G(t,\bU)=\bld{F}(\bU).
\)

An explicit RK discretization of
\eqref{IMEX-RK-CSF: unified ODE form} is straightforward, but the stiffness
of the semidiscrete geometric-flow system can impose severe time-step
restrictions. We therefore adopt a semi-implicit treatment: the nonlinear
matrix $\bld{F}(\bU)$ is evaluated at already available stage values,
whereas the stage solution on which it acts may be treated implicitly.
As a result, each stage requires only the solution of a linear system.

Motivated by the strategy in \cite{song2016energy}, we consider the following
$p$-stage IMEX--RK scheme:
\begin{align}
\begin{split}
\bU^{m,0}
&=
\bU^m,
\\
\bU^{m,i}
&=
\sum_{k=0}^{i-1}
\alpha_{i-1,k}\bU^{m,k}
+
\tau
\sum_{\ell=0}^{i}
\sum_{k=0}^{i-1}
\beta_{\ell k}^{(i)}
\bld{F}(\bU^{m,k})\bU^{m,\ell},
\qquad
i=1,\ldots,p,
\\
\bU^{m+1}
&=
\bU^{m,p}.
\end{split}
\label{IMEX-RK-CSF: Form of IMEX-RK scheme}
\end{align}
The coefficients are collected in
\[
C_p
=
\begin{pmatrix}
\alpha_{00} & 0 & \cdots & 0 \\
\alpha_{10} & \alpha_{11} & \cdots & 0 \\
\vdots & \vdots & \ddots & \vdots \\
\alpha_{p-1,0} & \alpha_{p-1,1} & \cdots & \alpha_{p-1,p-1}
\end{pmatrix},
\quad
B_p^{(i)}
=\begin{pmatrix}
\beta_{00}^{(i)}
& \beta_{01}^{(i)}
& \cdots
& \beta_{0,i-1}^{(i)}
\\
\beta_{10}^{(i)}
& \beta_{11}^{(i)}
& \cdots
& \beta_{1,i-1}^{(i)}
\\
\vdots
& \vdots
& \ddots
& \vdots
\\
\beta_{i0}^{(i)}
& \beta_{i1}^{(i)}
& \cdots
& \beta_{i,i-1}^{(i)}
\end{pmatrix},
\]
for $i=1,\ldots,p$.
Consistency at the first stage requires $\alpha_{00}=1$.

\begin{remark}
It was shown in \cite{song2016energy} that no two-stage, second-order method
exists within the particular framework considered there. Our formulation
introduces the additional coefficients contained in $B_p^{(i)}$, which
provides greater flexibility in solving the order conditions and allows us
to construct a two-stage, second-order scheme.
\end{remark}

\medskip
\begin{remark}		
Many standard IMEX--RK methods are designed for an additive splitting of the
form
\[
\bU_t=\mathcal{L}\bU+\bld{N}(\bU),
\]
where $\mathcal{L}$ is linear and the nonlinear term $\bld{N}(\bU)$ is
treated explicitly. Such a splitting does not directly reflect the
structure of \eqref{IMEX-RK-CSF: unified ODE form}, because
$\bld{F}(\bU)$ itself depends nonlinearly on $\bU$.
Ding \cite{ding2025semi} considered systems of the more general form
\[
\bU_t
=
\bld{f}(t,\bU)+G(t,\bU)\bU
\]
and introduced a semi-implicit treatment at intermediate stages.
In the present construction, every stage, including the final one, retains
this semi-implicit structure.
\end{remark}

To study linear stability, we apply
\eqref{IMEX-RK-CSF: Form of IMEX-RK scheme} to the test equation
\[
u_t=\lambda u,
\qquad
\Re(\lambda)\leq0.
\]
Writing $z=\lambda\tau$, the resulting update has the form
\(
u^{m+1}=R(z)u^m,
\)
where $R$ is the stability function. We seek coefficient sets for which the
method is $A$-stable, namely,
\[
|R(z)|<1
\qquad
\text{for all }z\in\mathbb{C}
\text{ with }\Re(z)<0.
\]
In particular, $R$ must have no nonremovable poles in the open left
half-plane. The order conditions are obtained by matching Taylor expansions.
We next construct schemes of second and third order.

\subsection{Two-stage second-order IMEX-RK methods}
\
We first set $p=2$ and require second-order accuracy. Taylor expansion gives
the following order conditions:
\begin{align*}
\alpha_{00}
&=1,
\\
\alpha_{10}+\alpha_{11}
&=1,
\\
\beta_{20}^{(2)}
+\beta_{21}^{(2)}
+\alpha_{11}
\left(
\beta_{00}^{(1)}+\beta_{10}^{(1)}
\right)
+\beta_{00}^{(2)}
+\beta_{01}^{(2)}
+\beta_{10}^{(2)}
+\beta_{11}^{(2)}
&=1,
\\
\left(
\beta_{21}^{(2)}
+\beta_{01}^{(2)}
+\beta_{11}^{(2)}
\right)
\left(
\beta_{00}^{(1)}
+\beta_{10}^{(1)}
\right)
&=\frac12,
\\
\beta_{20}^{(2)}
+\beta_{21}^{(2)}
+
\left(
\beta_{00}^{(1)}
+\beta_{10}^{(1)}
\right)
\left(
\alpha_{11}\beta_{10}^{(1)}
+\beta_{10}^{(2)}
+\beta_{11}^{(2)}
\right)
&=\frac12.
\end{align*}
For example, five coefficient sets satisfying the second-order conditions are listed in
Table~\ref{IMEX-RK: coefficients 2-stage}.

\begin{table*}[!h]
	\centering \begin{tabular*}{\textwidth}{@{\extracolsep{\fill}}ccccccccccc}
		\toprule
		& $\alpha_{10}$ & $\alpha_{11}$ & $\beta_{00}^{(1)}$ & $\beta_{10}^{(1)}$ & $\beta_{00}^{(2)}$ & $\beta_{01}^{(2)}$ & $\beta_{10}^{(2)}$  & $\beta_{11}^{(2)}$ & $\beta_{20}^{(2)}$ & $\beta_{21}^{(2)}$\\
		\hline
		(1) & 0 & 1 & 1  & 1 & $1/4$ & $1/4$ & 0 & 0  & $-3/2$ & 0\\
		\hline
		(2) & 0 & 1 & 1 & 1 & $-1/2$ & $1/8$ & $-7/8$ & 0 & $1/8$ & $1/8$\\
		\hline
		(3) & 0  & 1 & 0 & 2 & $1/8$ & $ 1/8$ & $ -9/4$ & 0& $7/8$ & $1/8$\\
		\hline
		(4) & $1$ & $0$ & $0$ & $1/2$ & 0 & $1/2$ & 0 & 0 & 0 & $1/2 $
		\\
\hline
(5) & $1/2$ & $1/2$ & $2/3$ & $1/3$ & 0 & $1/6$ & 0 & $-1/6$ & 0 & $1/2 $
				\\
		\bottomrule
	\end{tabular*}
\caption{Five coefficient sets satisfying the second-order conditions.}
\label{IMEX-RK: coefficients 2-stage}
\end{table*}

Tedious calculations give the stability function
\begin{align*}
&R(z)=
\frac{1+P_1z+P_2z^2}
{
\left(1-\beta_{10}^{(1)}z\right)
\left[
1-
\left(
\beta_{20}^{(2)}+\beta_{21}^{(2)}
\right)z
\right]
},
\\
&P_1 =
\beta_{00}^{(2)}
+\beta_{01}^{(2)}
+\beta_{10}^{(2)}
+\beta_{11}^{(2)}
+
\alpha_{11}
\left(
\beta_{00}^{(1)}+\beta_{10}^{(1)}
\right)
-
\beta_{10}^{(1)},
\\
&P_2
=
-
\left(
\beta_{00}^{(2)}+\beta_{01}^{(2)}
\right)
\beta_{10}^{(1)}
+
\left(
\beta_{10}^{(2)}+\beta_{11}^{(2)}
\right)
\beta_{00}^{(1)}.
\end{align*}
It can be easily checked that coefficient sets (1) and (2) do not satisfy the $A$-stability criterion, whereas sets (3)-(5) are $A$-stable.

Applying the 2-stage method to
\eqref{IMEX-RK-CSF: unified ODE form} gives
\begin{equation}\label{RKG}
\begin{split}
\bU^{m,1}
={}&
\bU^m
-
\beta_{00}^{(1)}
\tau
A(\bU^m)^{-1}M(\bU^m)\bU^m
-
\beta_{10}^{(1)}
\tau
A(\bU^m)^{-1}M(\bU^m)\bU^{m,1},
\\
\bU^{m+1}
={}&
\alpha_{10}\bU^m
+
\alpha_{11}\bU^{m,1}
\\
&-
\tau
\left[
\beta_{00}^{(2)}
A(\bU^m)^{-1}M(\bU^m)
+
\beta_{01}^{(2)}
A(\bU^{m,1})^{-1}M(\bU^{m,1})
\right]\bU^m
\\
&-
\tau
\left[
\beta_{10}^{(2)}
A(\bU^m)^{-1}M(\bU^m)
+
\beta_{11}^{(2)}
A(\bU^{m,1})^{-1}M(\bU^{m,1})
\right]\bU^{m,1}
\\
&-
\tau
\left[
\beta_{20}^{(2)}
A(\bU^m)^{-1}M(\bU^m)
+
\beta_{21}^{(2)}
A(\bU^{m,1})^{-1}M(\bU^{m,1})
\right]\bU^{m+1}.
\end{split}
\end{equation}

For Dziuk's formulation, the first stage is equivalently written as
\begin{align}
\left(
\frac{\bX^{m,1}-\bX^m}{\tau},
\bphi_h
\right)^h_{\Gamma^m}
+
\left(
\p_s
\left(
\beta_{00}^{(1)}\bX^m
+
\beta_{10}^{(1)}\bX^{m,1}
\right),
\p_s\bphi_h
\right)_{\Gamma^m}=0.
\label{IMEX-RK-CSF: Dziuk 2stage eq 1}
\end{align}
For the stage-local formulation used later in the BGN extension, we impose
the reduced-coefficient conditions
\[
\beta_{00}^{(2)}
=
\beta_{10}^{(2)}
=
\beta_{20}^{(2)}
=
0.
\]
Coefficient sets (4)-(5) in
Table~\ref{IMEX-RK: coefficients 2-stage} satisfies these conditions.
The second stage then reduces to
\begin{equation}\label{IMEX-RK-CSF: Dziuk 2stage eq 2}
	\frac{1}{\tau}\bigg( \bX^{m+1} - \sum\limits_{j=0}^{1}\alpha_{1j}\bX^{m,j},\bphi_h\bigg)^h_{\Gamma^{m,1}} + \bigg( \p_s\bigg( \sum_{k=0}^2\beta_{k1}^{(2)}\bX^{m,k} \bigg), \p_s\bphi_h\bigg)_{\Gamma^{m,1}} = 0.
\end{equation}

\begin{remark}
The reduced-coefficient structure allows the general scheme \eqref{RKG} to be reformulated as a parametric finite element scheme of the form \eqref{IMEX-RK-CSF: Dziuk 2stage eq 2}, in which both the variational formulation and the associated integrals are posed on the curve from the immediately preceding stage. This feature is crucial for directly extending IMEX--RK methods to BGN-type formulations in the next section.
Particularly, coefficient set (4) in
Table~\ref{IMEX-RK: coefficients 2-stage} reduces
\eqref{IMEX-RK-CSF: s-stage CSF+DD} to the predictor--corrector scheme
of \cite{deckelnick2026second}.
\end{remark}

\subsection{Three-stage third-order IMEX-RK methods}
\
We next set $p=3$ and construct a third-order scheme.
The full third-order conditions obtained by Taylor expansion are lengthy.
To obtain a stage-local formulation, we impose
\[
\beta_{00}^{(2)}
=
\beta_{10}^{(2)}
=
\beta_{20}^{(2)}
=
0,\quad
\beta_{00}^{(3)}
=
\beta_{01}^{(3)}
=
\beta_{10}^{(3)}
=
\beta_{11}^{(3)}
=
\beta_{20}^{(3)}
=
\beta_{21}^{(3)}
=
0.
\]
Under these restrictions, the 3-stage method applied to
\eqref{IMEX-RK-CSF: unified ODE form} becomes
\begin{align*}
	\bU^{m,1} &= \bU^m - \beta_{00}^{(1)} \tau A(\bU^m)^{-1} M(\bU^m) \bU^m - \beta_{10}^{(1)} \tau  A(\bU^m)^{-1} M(\bU^m) \bU^{m,1}, \\
	\bU^{m,2} &= \alpha_{10}\bU^{m} + \alpha_{11}\bU^{m,1} - \tau \beta_{01}^{(2)} A(\bU^{m,1})^{-1} M(\bU^{m,1}) \bU^m \\
	&\quad - \tau A(\bU^{m,1})^{-1} M(\bU^{m,1}) \left( \beta_{11}^{(2)} \bU^{m,1} +\beta_{21}^{(2)} \bU ^{m,2} \right),\\
	\bU^{m+1} &= \bU^{m+1,3}=\alpha_{20}\bU^{m} + \alpha_{21}\bU^{m,1} + \alpha_{22} \bU^{m,2} \\
	&\quad - \tau A(\bU^{m,2})^{-1} M(\bU^{m,2}) \left( \beta_{02}^{(3)}  \bU^m  + \beta_{12}^{(3)} \bU^{m,1} +\beta_{22}^{(3)} \bU ^{m,2} +\beta_{32}^{(3)} \bU ^{m+1}\right).
\end{align*}
One coefficient set satisfying the third-order conditions is
\begin{align}
	\begin{split}
		\alpha_{10} &= 0.75025649618, \quad
		\alpha_{11} = 0.24974350382, \\
		\alpha_{20} &= 0.35939205428, \quad
		\alpha_{21} = 0.63891460578, \quad
		\alpha_{22} = 0.00169333994, \\
		\beta_{10}^{(1)}  &= 0.39077142550, \quad
		\beta_{11}^{(2)}  = 0.03044132237, \quad
		\beta_{21}^{(2)}  = 0.53884073494, \\
		\beta_{12}^{(3)}  &= 0.83687560240, \quad
		\beta_{22}^{(3)}  = -0.48747543191, \quad
		\beta_{32}^{(3)}  = 0.3998010127,
	\end{split} \label{IMEX-RK: Coefficients 3-stage}
\end{align}
with $\beta_{00}^{(1)} = \beta_{01}^{(2)} = \beta_{02}^{(3)} = 0$. The coefficients are displayed to 11 significant digits. It can be checked that the resulting scheme is $A$-stable.

More generally, under the stage-local coefficient structure, the
$p$-stage schemes for Dziuk's and Deckelnick--Dziuk's formulations can be
written, respectively, as
\begin{align}
	\left( D_i\bX^m,\bphi_h\right)^h_{\Gamma^{m,i-1}} + \left( \p_s \bY^{m,i}, \p_s\bphi_h\right)_{\Gamma^{m,i-1}} &= 0,\label{IMEX-RK-CSF: s-stage CSF+Dziuk}\\
	\left( |\bX^{m,i-1}_\rho|^2 D_i\bX^m,\bphi_h\right)^h_{\S^1} + \left( \p_\rho \bY^{m,i}, \p_\rho\bphi_h \right)_{\S^1} &= 0, \label{IMEX-RK-CSF: s-stage CSF+DD}
\end{align}
for every $\bphi_h\in[V^h]^2$ and $i=1,\ldots,p$,
where  \begin{equation}\label{stageq}
D_i\bX^m
:=
\frac{1}{\tau}\left(
\bX^{m,i}
-\displaystyle\sum_{k=0}^{i-1}
\alpha_{i-1,k}\bX^{m,k}\right),
\qquad
\bY^{m,i}
:=
\sum_{k=0}^{i}
\beta_{k,i-1}^{(i)}\bX^{m,k}.
\end{equation}
with $\bX^{m,0}=\bX^m$, $\Gamma^{m,i}=\bX^{m,i}(\S^1)$.

We conclude this section with a well-posedness result for the stage equations.

\btheorem[Well-posedness]\label{IMEX-RK-CSF: Wellposedness of Dziuk-DD+CSF} Suppose that, for each $i=1,\ldots,p$,
$\beta_{i,i-1}^{(i)}\geq0$,
and the preceding-stage polygon is nondegenerate:
$\min\limits_{1\leq j\leq N}
|\bld{h}_j^{m,i-1}|
>0,$ where
$\bld{h}_j^{m,i-1}
=\bX^{m,i-1}(\rho_j)
-\bX^{m,i-1}(\rho_{j-1}).$
Then every stage of
\eqref{IMEX-RK-CSF: s-stage CSF+Dziuk}
admits a unique solution in $[V^h]^2$. In particular,
$\bX^{m+1}=\bX^{m,p}$ is uniquely determined.
The same conclusion holds for
\eqref{IMEX-RK-CSF: s-stage CSF+DD}.
\etheorem
\bproof

We argue stage by stage. Suppose that
$\bX^{m,0},\ldots,\bX^{m,i-1}$ are already known.
The $i$th stage is a finite-dimensional linear system. To prove uniqueness,
it suffices to show that its associated homogeneous problem has only the
zero solution. For Dziuk's formulation, the homogeneous problem is
\[\frac{1}{\tau}
\left(\bX,
\bphi_h
\right)^h_{\Gamma^{m,i-1}}
+
\beta_{i,i-1}^{(i)}
\left(
\p_s\bX,
\p_s\bphi_h
\right)_{\Gamma^{m,i-1}}
=
0,
\qquad
\forall\,\bphi_h\in[V^h]^2.
\]
Choosing $\bphi_h=\bX$ and using the nondegeneracy assumption together with $\beta_{i,i-1}^{(i)}\geq 0$, we obtain $\bX=0$. Hence the stage matrix is nonsingular, and the $i$th stage has a unique
solution. Induction over $i=1,\ldots,p$ proves the result for
\eqref{IMEX-RK-CSF: s-stage CSF+Dziuk}.

For \eqref{IMEX-RK-CSF: s-stage CSF+DD}, the same argument applies because
$|\p_\rho\bX^{m,i-1}|^2>0$ elementwise under the nondegeneracy assumption.
\eproof

\begin{remark}
	The coefficient set \eqref{IMEX-RK: Coefficients 3-stage} satisfies
$\beta_{i,i-1}^{(i)}>0$, $i=1,2,3$, and hence satisfies the coefficient condition in the preceding theorem.
\end{remark}

\section{Extension to BGN methods}\label{sec: To BGN}
\
We now extend the IMEX--RK framework developed in the previous section to BGN-type formulations. The reduced-coefficient structure is essential for
this extension: at the $i$th RK stage, all variational forms and geometric quantities are evaluated on the preceding stage $\Gamma^{m,i-1}$. Thus, $\Gamma^{m,i-1}$ naturally serves as the prediction
curve for the $i$th stage. We first present the construction for CSF and then for CDF.

\subsection{Curve shortening flow}
\
The BGN formulation of CSF is based on
\begin{align*}
\p_t\bX\cdot\bld{N}
=-\kappa,
\qquad
\kappa\bld{N}
=-\p_{ss}\bX.
\end{align*}
The corresponding weak
formulation is
\begin{align*}
\left(
\p_t\bX\cdot\bld{N},
\varphi
\right)_{\Gamma(t)}
+
\left(
\kappa,\varphi
\right)_{\Gamma(t)}
&=
0,
\qquad
\forall\,\varphi\in H^1(\Gamma(t)),
\\
\left(
\kappa\bld{N},
\bpsi
\right)_{\Gamma(t)}
-
\left(
\p_s\bX,
\p_s\bpsi
\right)_{\Gamma(t)}
&=
0,
\qquad
\forall\,\bpsi\in[H^1(\Gamma(t))]^2.
\end{align*}

For $\bX_h\in[V^h]^2$, recall that
$\bld{N}_h|_{I_j}
=-\frac{\bld{h}_j^\perp}{|\bld{h}_j|}$.
The semidiscrete BGN scheme for CSF is then formulated as follows.
Given $\bX_h(\cdot,0)\in[V^h]^2$, find
\[
\bX_h(\cdot,t)\in[V^h]^2,
\qquad
\kappa_h(\cdot,t)\in V^h,
\]
such that
\begin{subequations}\label{IMEX-RK-CSF: BGN semi-discrete}
	\begin{align}
		\begin{split}
			\lb \p_t\bX_h\cdot\bld{N}_h,\varphi_h\rb^h_{\Gamma_h(t)} + \lb \kappa_h,\varphi_h\rb^h_{\Gamma_h(t)} & = 0, \quad \forall \,  \varphi_h \in V^h ,\label{IMEX-RK-CSF: BGN semi-discrete eq1}
		\end{split} \\
		\begin{split}
			\lb\kappa_h\bld{N}_h,\bpsi_h\rb^h_{\Gamma_h(t)} - \lb \p_s\bX_h,\p_s\bpsi_h\rb_{\Gamma_h(t)} & = 0,\quad \forall \,  \bpsi_h\in [V^h]^2. \label{IMEX-RK-CSF: BGN semi-discrete eq2}
		\end{split}
	\end{align}
\end{subequations}
A semi-implicit Euler discretization gives the classical first-order BGN
scheme. Given $\bX^m\in[V^h]^2$, find
$\bX^{m+1}\in[V^h]^2$ and $\kappa^{m+1}\in V^h$ such that
\begin{subequations}
\label{IMEX-RK-CSF: BGN1 method}
\begin{align}
\frac{1}{\tau}\left(
(\bX^{m+1}-\bX^m)
\cdot\bld{N}^m,
\varphi_h
\right)^h_{\Gamma^m}
+
\left(
\kappa^{m+1},
\varphi_h
\right)^h_{\Gamma^m}
&=
0,
\label{IMEX-RK-CSF: BGN1 method eq1}
\\
\left(
\kappa^{m+1}\bld{N}^m,
\bpsi_h
\right)^h_{\Gamma^m}
-
\left(
\p_s\bX^{m+1},
\p_s\bpsi_h
\right)_{\Gamma^m}
&=
0,
\label{IMEX-RK-CSF: BGN1 method eq2}
\end{align}
\end{subequations}
for all $\varphi_h\in V^h$ and $\bpsi_h\in[V^h]^2$.

We next apply the stage-local IMEX--RK construction from the previous
section.
At the $i$th stage, the geometry is evaluated entirely on
$\Gamma^{m,i-1}$. We seek
\[
\bX^{m,i}\in[V^h]^2,
\qquad
\kappa^{m,i}\in V^h,
\]
such that
\begin{subequations}
\label{IMEX-RK-CSF: s-stage BGN+CSF}
\begin{align}
\left(
D_i\bX^m\cdot\bld{N}^{m,i-1},
\varphi_h
\right)^h_{\Gamma^{m,i-1}}
+
\left(
\kappa^{m,i},
\varphi_h
\right)^h_{\Gamma^{m,i-1}}
&=
0,
\label{IMEX-RK-CSF: s-stage BGN+CSF eq1}
\\
\left(
\kappa^{m,i}\bld{N}^{m,i-1},
\bpsi_h
\right)^h_{\Gamma^{m,i-1}}
-
\left(
\p_s\bY^{m,i},
\p_s\bpsi_h
\right)_{\Gamma^{m,i-1}}
&=
0,
\label{IMEX-RK-CSF: s-stage BGN+CSF eq2}
\end{align}
\end{subequations}
for all $\varphi_h\in V^h$, $\bpsi_h\in[V^h]^2$, and
$i=1,\ldots,p$. Finally, we set
$\bX^{m+1}=\bX^{m,p}$.

The structure of
\eqref{IMEX-RK-CSF: s-stage BGN+CSF}
is particularly simple. The curve $\Gamma^{m,i-1}$ obtained at the
preceding stage supplies all geometric quantities required at the current
stage, including the normal vector, the arc-length derivative, and the
integration measure. Hence, no separate prediction curve is required. Particularly, for $p=2$, coefficient set (4) in Table~\ref{IMEX-RK: coefficients 2-stage} reduces \eqref{IMEX-RK-CSF: s-stage BGN+CSF}
to the BGN predictor--corrector scheme in
\cite{jiang2025predictor}.

\btheorem[Well-posedness]\label{IMEX-RK-CSF: Wellposedness of BGN+CSF} For $m \geq 0$, we assume that  $\beta_{i,i-1}^{(i)} > 0$,  for each $i=1, \ldots, p$, and that the preceding-stage polygon $\Gamma^{m,i-1}$ satisfies
\begin{equation}\label{normalv}
\min_{1\leq j\leq N}
|\bld{h}_j^{m,i-1}|
>0,
\qquad
\dim
\operatorname{span}
\left\{
\bomega_j^{m,i-1}
\right\}_{j=1}^{N}
=2, \end{equation}
where the weighted vertex normal vector is defined by
\[
\bomega_j^{m,i}
:=-\frac{
(\bX^{m,i}(\rho_{j+1})-\bX^{m,i}(\rho_{j-1}))^\perp}{
|\bld{h}_j^{m,i}|+|\bld{h}_{j+1}^{m,i}|
},
\]
with the indices understood periodically.
Then every stage of
\eqref{IMEX-RK-CSF: s-stage BGN+CSF}
has a unique solution
$(\bX^{m,i},\kappa^{m,i})
\in[V^h]^2\times V^h$. Consequently,
$\bX^{m+1}=\bX^{m,p}$ is uniquely determined.
\etheorem

\bproof
The proof follows the same argument in \cite[Theorem 76]{barrett2020parametric} for classical BGN method and we omit it here for brevity.
\eproof

\subsection{Curve diffusion flow}
\
We next consider curve diffusion flow (CDF), whose BGN formulation reads as
\[
\p_t\bX\cdot\bld{N}
=\p_{ss}\kappa,\qquad
\kappa\bld{N}=-\p_{ss}\bX.
\]
The weak formulation leads to the BGN
semidiscrete scheme: find
$\bX_h(\cdot,t)\in[V^h]^2$ and
$\kappa_h(\cdot,t)\in V^h$ such that
\begin{subequations}\label{IMEX-RK-CDF: BGN semi-discrete}
	\begin{align}
		\begin{split}
			\lb \p_t\bX_h\cdot\bld{N}_h,\varphi_h\rb^h_{\Gamma_h(t)} + \lb\p_s\kappa_h,\p_s\varphi_h\rb^h_{\Gamma_h(t)} & = 0, \quad \forall \,  \varphi_h \in V^h ,\label{IMEX-RK-CDF: BGN semi-discrete eq1}
		\end{split} \\
		\begin{split}
			\lb\kappa_h\bld{N}_h,\bpsi_h\rb^h_{\Gamma_h(t)} - \lb \p_s\bX_h,\p_s\bpsi_h\rb_{\Gamma_h(t)} & = 0,\quad \,  \forall \bpsi_h\in [V^h]^2, \label{IMEX-RK-CDF: BGN semi-discrete eq2}
		\end{split}
	\end{align}
\end{subequations}
The classical first-order BGN discretization is obtained by finding $\bX^{m+1}\in[V^h]^2$ and $\kappa^{m+1}\in V^h$ such that
\begin{align}
    \begin{split}
       \frac{1}{\tau} \lb \bX^{m+1} - \bX^m,\varphi_h\bld{N}^m\rb^h_{\Gamma^m} + \lb \p_s\kappa^{m+1},\p_s\varphi_h\rb_{\Gamma^m} & = 0,\quad \forall\, \varphi_h \in V^h,\\
        \lb\kappa^{m+1}\bld{N}^m,\bpsi_h\rb^h_{\Gamma^m} -\lb\p_s\bX^{m+1},\p_s\bpsi_h\rb_{\Gamma^m} & = 0,\quad \forall \, \bpsi_h \lmb V^h\rmb^2.
    \end{split}
\end{align}
Using the same stage quantities $D_i\bX^m$ and $\bY^{m,i}$ introduced
above, the $p$-stage IMEX--RK discretization is
\begin{subequations}
\label{IMEX-RK-CDF: s-stage BGN+CDF}
\begin{align}
\left(
D_i\bX^m\cdot\bld{N}^{m,i-1},
\varphi_h
\right)^h_{\Gamma^{m,i-1}}
+
\left(
\p_s\kappa^{m,i},
\p_s\varphi_h
\right)_{\Gamma^{m,i-1}}
&=
0,
\label{IMEX-RK-CDF: s-stage BGN+CDF eq1}
\\
\left(
\kappa^{m,i}\bld{N}^{m,i-1},
\bpsi_h
\right)^h_{\Gamma^{m,i-1}}
-
\left(
\p_s\bY^{m,i},
\p_s\bpsi_h
\right)_{\Gamma^{m,i-1}}
&=
0,
\label{IMEX-RK-CDF: s-stage BGN+CDF eq2}
\end{align}
\end{subequations}
for all $\varphi_h\in V^h$, $\bpsi_h\in[V^h]^2$, and
$i=1,\ldots,p$.

Under the assumptions of
Theorem~\ref{IMEX-RK-CSF: Wellposedness of BGN+CSF},
the stage equations
\eqref{IMEX-RK-CDF: s-stage BGN+CDF}
are also uniquely solvable. Indeed, the corresponding homogeneous problem
yields
\[
\tau
\|\p_s\kappa\|_{\Gamma^{m,i-1}}^2
+
\beta_{i,i-1}^{(i)}
\|\p_s\bX\|_{\Gamma^{m,i-1}}^2
=
0.
\]
Thus both $\kappa$ and $\bX$ are constant. The two homogeneous equations,
together with the vertex-normal condition in
\eqref{normalv}, then imply that both constants
vanish.

\begin{remark}
A central issue in high-order BGN discretizations is the choice of the
geometry on which each time-discrete variational problem is posed.
For the BDF$k$ schemes in \cite{jiang2024stable}, sufficiently accurate
prediction curves are required to retain high-order consistency and robust
mesh redistribution. In the present IMEX--RK formulation, this prediction is
provided intrinsically by the multistage structure: the curve
$\Gamma^{m,i-1}$ generated at the preceding stage serves directly as the
integration and prediction curve for stage $i$. Consequently, no separate
prediction procedure or additional refined time stepping is required.
\end{remark}

\section{Extension to surface case}
We now extend the proposed IMEX--RK framework to
two-dimensional surfaces in $\mathbb{R}^3$. The construction follows the same stage-local principle developed in the previous sections. At each RK
stage, all geometric quantities are evaluated on the surface generated by the immediately preceding stage. We consider mean curvature flow (MCF) and surface
diffusion flow (SDF).

\subsection{Dziuk's method for mean curvature flow}

Let $\Sigma$ be a smooth closed reference surface, and let
$\bX(\cdot,t):\Sigma\rightarrow\mathbb{R}^3$
be a smooth immersion. We denote the evolving surface by $\Gamma(t)=\bX(\Sigma,t)$.
Let $\bld{N}$ be its outward unit normal and
$\Delta_{\Gamma}$ the Laplace--Beltrami operator. With the convention
$\mH
=-\Delta_{\Gamma}\bX\cdot\bld{N}$,
the mean curvature vector satisfies
$\Delta_{\Gamma}\bX=-\mH\bld{N}$.
The $L^2$-gradient flow of the surface area
$\mathcal{A}(\Gamma)
=\int_{\Gamma}1\,\upd\sigma$
is therefore
\begin{align}
\p_t\bX
=
-\mH\bld{N}
=
\Delta_{\Gamma}\bX.
\label{IMEX-RK-MCF: MCF in geometric form}
\end{align}

We next recall the spatial finite element setting. Let $\Gamma_h(t)=\bigcup_{j=1}^{J}\overline{\sigma_j(t)}$ be a
conforming polyhedral approximation of $\Gamma(t)$,
where $\{\sigma_j(t)\}_{j=1}^{J}$ is a triangulation consisting of
nondegenerate triangles. We set
\[h(t)
=\max\limits_{1\leq j\leq J}
\operatorname{diam}\bigl(\sigma_j(t)\bigr).\]
We suppress the dependence on $t$ whenever no confusion can arise.

The continuous piecewise linear finite element space on $\Gamma_h$ is
\[
W^h(\Gamma_h)
:=
\left\{
u\in C(\Gamma_h):
u|_{\sigma_j}\ \text{is linear},
\quad
j=1,\ldots,J
\right\}.
\]
For piecewise continuous scalar- or vector-valued functions $u$ and $v$, we
define the $L^2$ inner product $(\cdot, \cdot)_{\Gamma_h}$ and its discrete version as
\[
(u,v)_{\Gamma_h}
=\int_{\Gamma_h}u\cdot v\,\upd\sigma,\qquad
(u,v)^h_{\Gamma_h}
:=\frac13
\sum_{j=1}^{J}
|\sigma_j|
\sum_{k=1}^{3}
\left.
(u\cdot v)
\right|_{\sigma_j}
(\bld{q}_{j_k}),
\]
where
$\bld{q}_{j_1}$, $\bld{q}_{j_2}$, and $\bld{q}_{j_3}$
are the vertices of $\sigma_j$.

We orient each triangle consistently with the outward normal and define
\[
\bld{N}_h|_{\sigma_j}
=
\bld{N}_j
:=
\frac{
(\bld{q}_{j_2}-\bld{q}_{j_1})
\times
(\bld{q}_{j_3}-\bld{q}_{j_1})
}{
\left|
(\bld{q}_{j_2}-\bld{q}_{j_1})
\times
(\bld{q}_{j_3}-\bld{q}_{j_1})
\right|
}.
\]

Dziuk's semidiscrete PFEM for MCF takes the form
\begin{align}
\left(
\p_t\bX_h,
\bphi_h
\right)^h_{\Gamma_h(t)}
+
\left(
\nabla_{\Gamma_h}\bX_h,
\nabla_{\Gamma_h}\bphi_h
\right)_{\Gamma_h(t)}
=
0,
\qquad
\forall\,\bphi_h\in[W^h(\Gamma_h(t))]^3,
\label{IMEX-RK-MCF: Dziuk's scheme}
\end{align}
where $\nabla_\Gamma$ is the surface gradient piecewisely defined over the polyhedra $\Gamma$.
 Applying the $p$-stage IMEX--RK discretization to
\eqref{IMEX-RK-MCF: Dziuk's scheme} gives
\begin{align}
	\left(D_i \bX^m,\bphi_h\right)^h_{\Gamma^{m,i-1}} + \left( \nabla_\Gamma \bY^{m,i}, \nabla_\Gamma\bphi_h\right)_{\Gamma^{m,i-1}} = 0,\label{IMEX-RK-MCF: s-stage MCF+Dziuk}
\end{align}
for all $\bphi_h\in [W^h(\Gamma^{m,i-1})]^3$ for each $i$-th stage, where $i=1,\ldots,p$.

\subsection{BGN methods for mean curvature flow and surface diffusion flow}
\
The BGN formulation of MCF separates the normal velocity from the
mean-curvature relation:
\[
\p_t\bX\cdot\bld{N}
=-\mH,
\qquad
\mH\bld{N}=
-\Delta_{\Gamma}\bX.
\]
Using the same stage-local construction as for CSF, the $p$-stage BGN
discretization of MCF is: for $i=1,\ldots,p$, find
\(
\bX^{m,i}\in[W^h(\Gamma^{m,i-1})]^3\),
\(\mH^{m,i}\in W^h(\Gamma^{m,i-1}), \)
such that
\begin{subequations}
\label{IMEX-RK-MCF: s-stage BGN+MCF}
\begin{align}
\left(
D_i\bX^m\cdot\bld{N}^{m,i-1},
\varphi_h
\right)^h_{\Gamma^{m,i-1}}
+
\left(
\mH^{m,i},
\varphi_h
\right)^h_{\Gamma^{m,i-1}}
&=
0,
\label{IMEX-RK-MCF: s-stage BGN+MCF eq1}
\\
\left(
\mH^{m,i}\bld{N}^{m,i-1},
\bpsi_h
\right)^h_{\Gamma^{m,i-1}}
-
\left(
\nabla_{\Gamma^{m,i-1}}\bY^{m,i},
\nabla_{\Gamma^{m,i-1}}\bpsi_h
\right)_{\Gamma^{m,i-1}}
&=0,
\label{IMEX-RK-MCF: s-stage BGN+MCF eq2}
\end{align}
\end{subequations}
for all $\varphi_h\in W^h(\Gamma^{m,i-1})$,
$\bpsi_h\in[W^h(\Gamma^{m,i-1})]^3$.

We finally consider surface diffusion flow,
\begin{align}
\p_t\bX\cdot\bld{N}
=\Delta_{\Gamma}\mH,\qquad \mH\bld{N}
=-\Delta_{\Gamma}\bX.
\label{IMEX-RK-SDF: SDF}
\end{align}
The corresponding $p$-stage IMEX--RK scheme is: for
$i=1,\ldots,p$, find
$\bX^{m,i}\in[W^h(\Gamma^{m,i-1})]^3$,
$\mH^{m,i}\in W^h(\Gamma^{m,i-1})$, such that
\begin{subequations}
\label{IMEX-RK-SDF: s-stage BGN+SDF}
\begin{align}
\left(
D_i\bX^m\cdot\bld{N}^{m,i-1},
\varphi_h
\right)^h_{\Gamma^{m,i-1}}
+
\left(
\nabla_{\Gamma^{m,i-1}}\mH^{m,i},
\nabla_{\Gamma^{m,i-1}}\varphi_h
\right)_{\Gamma^{m,i-1}}
&=
0,
\label{IMEX-RK-SDF: s-stage BGN+SDF eq1}
\\
\left(
\mH^{m,i}\bld{N}^{m,i-1},
\bpsi_h
\right)^h_{\Gamma^{m,i-1}}
-
\left(
\nabla_{\Gamma^{m,i-1}}\bY^{m,i},
\nabla_{\Gamma^{m,i-1}}\bpsi_h
\right)_{\Gamma^{m,i-1}}
&=0,
\label{IMEX-RK-SDF: s-stage BGN+SDF eq2}
\end{align}
\end{subequations}
for all $\varphi_h\in W^h(\Gamma^{m,i-1})$, $\bpsi_h\in[W^h(\Gamma^{m,i-1})]^3$.

\section{Numerical results}
\

In this section, we assess the numerical performance of the proposed IMEX--RK schemes. The experiments serve four main purposes. First, we examine the role of $A$-stability in the robust simulation of geometric flows. Second, we verify the designed temporal orders of convergence for CSF, CDF and MCF. Third, we investigate whether the proposed high-order schemes retain the mesh-redistribution properties of the underlying parametric finite element formulations. Finally, we apply the methods to benchmark surface evolutions involving extinction and pinch-off singularities.

\begin{example}[$A$-stability]\label{example: A stablity}
	This example shows the importance of $A$-stability for robust geometric-flow simulations.
\end{example}

\begin{figure}[h!]
	\vspace{-2mm} \includegraphics[width=0.98\linewidth]{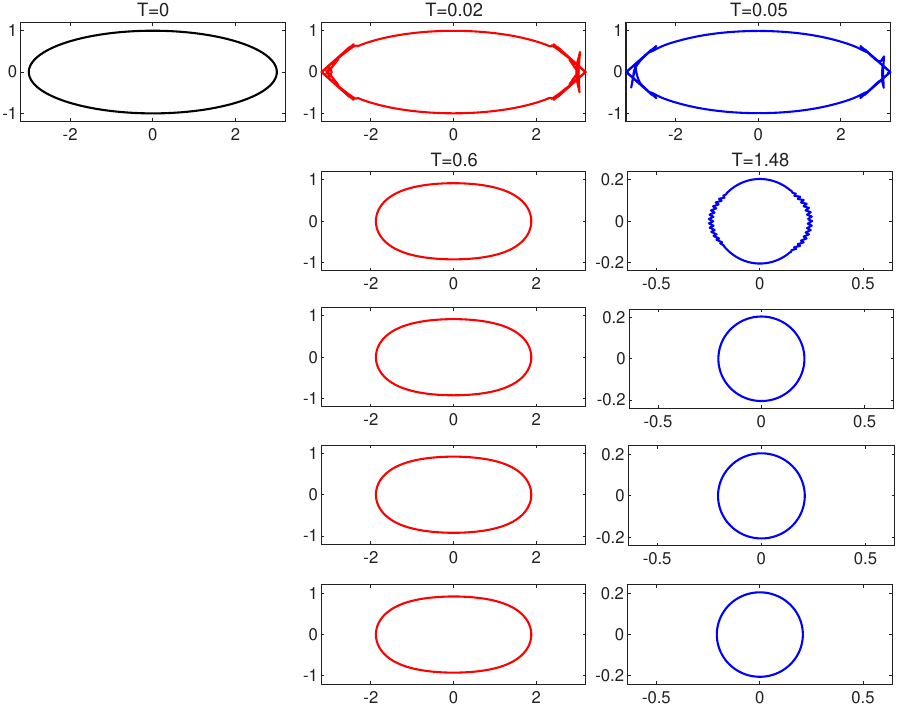}
	\vspace{-1mm}
	\caption{Evolution of the curve for the four coefficient sets. From top to bottom, the rows correspond to coefficient sets (1)-(5) in Table \ref{IMEX-RK: coefficients 2-stage}, respectively.}
	\label{fig: A stability}
\end{figure}

We examine the importance of $A$-stability of the methods by taking CSF as an example. The initial curve is the ellipse $x^2+9y^2=9$.
We apply the 2-stage IMEX--RK method to the Deckelnick--Dziuk formulation, i.e., \eqref{IMEX-RK-CSF: s-stage CSF+DD} with $p=2$. We use $N=128$ spatial nodes and the time step $\tau=10^{-3}$. Figure~\ref{fig: A stability} reveals a clear distinction between the coefficient sets. The solution obtained with set (1) becomes unstable rapidly. The contour produced by set (2) also deteriorates as the curve approaches extinction. In contrast, the two $A$-stable choices, sets (3)-(5), remain stable throughout the computation. These results indicate that $A$-stability is important for robust time integration in this stiff geometric-flow setting, even though all four schemes have the same formal order of accuracy.

For convergence tests involving tangential mesh motion, errors based directly on a particular parametrization are not appropriate. Both the Deckelnick--Dziuk and BGN formulations allow tangential redistribution of the mesh nodes. We therefore use the reparametrization-invariant manifold distance introduced in \cite{zhao2021energy} and also used in \cite{jiang2024second}. For two closed polygonal curves $\Gamma_1$ and $\Gamma_2$, it is defined by
\begin{equation}
\mathcal{M}\left(\Gamma_1,\Gamma_2\right):=
\left|\left(\Omega_1\setminus\Omega_2\right)
\cup
\left(\Omega_2\setminus\Omega_1\right)
\right|=
|\Omega_1|+|\Omega_2|-2|\Omega_1\cap\Omega_2|,
\label{IMEX-RK: manifold distance}
\end{equation}
where $\Omega_1$ and $\Omega_2$ are the bounded regions enclosed by $\Gamma_1$ and $\Gamma_2$, respectively. Here, $|\Omega|$ denotes the two-dimensional Lebesgue measure of $\Omega$.

For the convergence tests below, we report the results obtained using coefficient set (5) in Table~\ref{IMEX-RK: coefficients 2-stage} for $p=2$ and the coefficients in \eqref{IMEX-RK: Coefficients 3-stage} for $p=3$.

\begin{example}[Order of convergence for CSF] We examine the temporal convergence of
\eqref{IMEX-RK-CSF: s-stage CSF+Dziuk},
\eqref{IMEX-RK-CSF: s-stage CSF+DD}, and
\eqref{IMEX-RK-CSF: s-stage BGN+CSF}
for CSF.  We consider the unit circle as the initial curve
    \[
    \bX_0(\rho)=(\cos\rho,\sin\rho),
    \qquad \rho\in\S^1,
    \]
    for which the exact solution is
    \[
    \bX(t,\rho)
    =
    \sqrt{1-2t}\,
    \lb \cos\rho,\sin\rho\rb.
    \]
To make the spatial discretization error negligible
compared with the temporal error, we choose $N=20000$ and uniformly partition $\S^1$ with
\[h=2\pi/N, \qquad \rho_j=jh,\qquad j=1,\ldots,N,
\]
where the periodic identification is understood. The initial finite element
approximation $\bX^0\in[V^h]^2$ is defined by
\(\bX^0(\rho_j)
=(\cos\rho_j,\sin\rho_j)\),
$j=1,\ldots,N$.
For \(
\star\in\{\mathrm{Dziuk},\mathrm{DD},\mathrm{BGN}\},
\)
let
\(
\Gamma_{p,\star}^{T/\tau}
=
\bX_{p,\star}^{T/\tau}(\S^1)
\)
denote the numerical curve at time $T$ obtained with the corresponding $p$-stage method and time step $\tau$, and let \(
\Gamma(T)=\bX(T,\S^1)
\)
be the exact curve. We define
\[\mathcal{E}^{p}_{\star}(T,\tau)
=\mathcal{M} \lb
\Gamma_{p,\star}^{T/\tau},
\Gamma(T)
\rb,
\quad
\operatorname{Order}
=\frac{\log\left(
\mathcal{E}^{p}_{\star}(T,\tau_1)/
\mathcal{E}^{p}_{\star}(T,\tau_2)
\right)
}{
\log(\tau_1/\tau_2)
}.
\]
Here,
$\bX_{p,\star}^{T/\tau}$ denotes the numerical solution generated by
\eqref{IMEX-RK-CSF: s-stage CSF+Dziuk},
\eqref{IMEX-RK-CSF: s-stage CSF+DD}, or
\eqref{IMEX-RK-CSF: s-stage BGN+CSF}
when
$\star=\mathrm{Dziuk}$, $\mathrm{DD}$, or $\mathrm{BGN}$,
respectively.
\end{example}

\begin{table*}[h]
	\centering
	\tabcolsep0.3cm
	\begin{tabular*}{\textwidth}{@{\extracolsep{\fill}}cccccc}
		\toprule
		$\tau$ & $1/20$&  $1/40$ & $1/80$ & $1/160$ & $1/320$       \\
		\hline
		$\mathcal{E}^{2}_{\text{Dziuk}}$  &  8.097E-03  & 2.050E-03   & 5.062E-04 &  1.249-04 & 3.097E-05     \\
				$\text{Order}$   &  -- -- & 1.982 & 2.018 & 2.019 & 2.012  \\
				\hline
				$\mathcal{E}^{2}_{\text{DD}}$  &  8.097E-03  & 2.050E-03   & 5.062E-04 &  1.249-04 & 3.097E-05       \\
				$\text{Order}$   &  -- -- & 1.982 & 2.018 & 2.019 & 2.012  \\
				\hline
				$\mathcal{E}^{2}_{\text{BGN}}$  &  8.097E-03  & 2.050E-03   & 5.061E-04 &  1.249-04 & 3.091E-05    \\
				$\text{Order}$   &  -- -- & 1.982 & 2.018 & 2.019 & 2.015 \\
		\bottomrule
	\end{tabular*}
\caption{\label{table: Dziuk-DD-BGN-CSF-IMEX-2nd, circle} Errors and observed convergence orders at \(T=0.4\) for the 2-stage schemes
applied to the Dziuk, Deckelnick--Dziuk, and BGN formulations of CSF.}
\end{table*}

\begin{table*}[htbp!]
	\centering
	\tabcolsep0.3cm \begin{tabular*}{\textwidth}{@{\extracolsep{\fill}}cccccc}
		\toprule
		$\tau$ & $1/20$&  $1/40$ & $1/80$ & $1/160$ & $1/320$            \\
		\hline
		$\mathcal{E}^{3}_{\text{Dziuk}}$  &  2.772E-03  & 4.371E-04   & 6.138E-05 &  8.120E-06 & 1.037E-06     \\
		$\text{Order}$   &  -- --  & 2.665 & 2.832 & 2.918 & 2.969  \\
		\hline
		$\mathcal{E}^{3}_{\text{DD}}$  & 2.772E-03  & 4.371E-04   & 6.138E-05 &  8.120E-06 & 1.037E-06     \\
		$\text{Order}$   &  -- --  & 2.665 & 2.832 & 2.918 & 2.969  \\
		\hline
		$\mathcal{E}^{3}_{\text{BGN}}$  &   2.772E-03  & 4.371E-04   & 6.132E-05 &  8.058E-06 & 9.747E-07     \\
		$\text{Order}$   &  -- --  & 2.665 & 2.834 & 2.928 & 3.047  \\
		\bottomrule
	\end{tabular*}
	\caption{\label{table: Dziuk-DD-BGN-CSF-IMEX-3rd, circle} Errors and observed convergence orders at $T=0.4$ for the 3-stage schemes applied to the Dziuk, Deckelnick--Dziuk, and BGN formulations of CSF.}
\end{table*}

The results in
Tables~\ref{table: Dziuk-DD-BGN-CSF-IMEX-2nd, circle}
and~\ref{table: Dziuk-DD-BGN-CSF-IMEX-3rd, circle}
show that the observed rates approach $2$ for the 2-stage schemes and $3$
for the 3-stage schemes. This confirms the designed temporal orders of
convergence. The three spatial formulations produce nearly identical errors
in this radially symmetric test. Since the exact evolution remains circular,
no substantial tangential mesh redistribution is required. Consequently,
the differences among the tangential motions induced by the Dziuk,
Deckelnick--Dziuk, and BGN formulations have essentially no effect on the
geometric error measured by the reparametrization-invariant distance
$\mathcal{M}$.

\begin{example}[Order of convergence for MCF]
\label{example: convergence MCF sphere}

To complement the convergence tests for planar curves, we now examine the
temporal accuracy of the proposed IMEX--RK schemes for the evolution of a
two-dimensional surface in three-dimensional space. We consider
\eqref{IMEX-RK-MCF: s-stage MCF+Dziuk} and
\eqref{IMEX-RK-MCF: s-stage BGN+MCF}
with $p=2$ and $p=3$ for MCF starting from the unit sphere.
The exact solution remains spherical for $t<1/4$, with radius
$r(t)=\sqrt{1-4t}$. Following \cite{duan2024new}, we measure the error at the final time $T$ by
\[
\mathcal{E}_{\star}^{p,\infty}(\tau,T)
=\max_{1\leq j\leq K}
\left|
|\bld{q}^{T/\tau}_j|-r(T)
\right|,
\] where $p=2,3$ and
$\star\in\{\mathrm{Dziuk},\mathrm{BGN}\}$.
Here, $K$ denotes the number of vertices of the triangulated surface.
\end{example}

	\begin{table*}[!h]
		\centering
		\tabcolsep0.3cm
		\begin{tabular*}{\textwidth}{@{\extracolsep{\fill}}cccccc}
			\toprule
			$\tau$ &  ${1}/{40}$ & ${1}/{80}$ &  ${1}/{120}$ &  ${1}/{160}$       \\
			\hline
		$\mathcal{E}^{2,\infty}_{\text{Dziuk}}$  & 2.875E-03 & 7.312E-04 & 3.241E-04 &  1.820E-04    \\
				$\text{Order}$ & ---  & 1.975 & 2.006 & 2.005    \\
				\hline
				$\mathcal{E}^{2,\infty}_{\text{BGN}}$  & 2.860E-03 & 7.162E-04 & 3.091E-04 &  1.670E-04   \\
				$\text{Order}$   & --- & 1.998 & 2.073 &  2.140     \\
			\bottomrule
		\end{tabular*}
	\caption{\label{table: Dziuk-BGN-MCF-IMEX-2nd} Errors and observed convergence orders for the 2-stage IMEX--RK schemes applied to the Dziuk and BGN formulations of MCF for the shrinking-sphere
test.}
	\end{table*}
	
	\begin{table*}[h]
		\centering
		\tabcolsep0.3cm \begin{tabular*}{\textwidth}{@{\extracolsep{\fill}}cccccc}
			\toprule
			$\tau$ & ${1}/{40}$ & ${1}/{60}$ &  ${1}/{80}$ &  ${1}/{100}$      \\
			\hline
			$\mathcal{E}^{3,\infty}_{\text{Dziuk}}$ & 9.87E-04 & 3.43E-04 & 1.57E-04 & 8.53E-05   \\
			$\text{Order}$  & --- & 2.608 & 2.706 & 2.745   \\
			\hline
			$\mathcal{E}^{3,\infty}_{\text{BGN}}$ & 9.73E-04 & 3.28E-04 & 1.46E-04 & 7.04E-05   \\
			$\text{Order}$  & --- & 2.680 & 2.898 & 3.160   \\
			\bottomrule
		\end{tabular*}
	\caption{\label{table: Dziuk-BGN-MCF-IMEX-3rd} Errors and observed convergence orders for the 3-stage IMEX--RK schemes applied to the Dziuk and BGN formulations of MCF for the shrinking-sphere
test.}
	\end{table*}

We set $T=0.2$ and use a sufficiently fine spatial mesh so that the spatial
discretization error is negligible compared with the temporal error. The mesh
contains $J=588956$ triangular elements and $K=294480$ vertices.
Tables~\ref{table: Dziuk-BGN-MCF-IMEX-2nd}
and~\ref{table: Dziuk-BGN-MCF-IMEX-3rd}
report the errors and observed convergence orders for $p=2$ and $p=3$,
respectively. For the 2-stage schemes, the observed rates approach $2$ as
the time step is refined. The 3-stage schemes exhibit the corresponding
third-order behavior. These results show that the high-order temporal
accuracy observed for planar curve evolution is retained for surface evolution under MCF.

\begin{example}[Order of convergence for CDF]
We next examine the convergence of the proposed IMEX--RK discretization
\eqref{IMEX-RK-CDF: s-stage BGN+CDF} for CDF.
The initial curve is the ellipse
\(
\frac{x^2}{9}+\frac{y^2}{4}=1,
\)
and the final time is set to $T=0.2$.
Since the exact solution is unavailable, we estimate the convergence order
by comparing numerical solutions on successive space--time refinement levels.
 The error $\mathcal{E}_{\mathcal{M}}$ and the convergence order are given by
\[
\mathcal{E}_{\mathcal{M}}(T, \tau_1, \tau_2) := \mathcal{M}\left(\Gamma_{h_1}^{T/\tau_1}, \Gamma_{h_2}^{T/\tau_2}\right), \quad \text{Order} = \log\left(\frac{\mathcal{E}_{\mathcal{M}}(T, \tau_1, \tau_2)}{\mathcal{E}_{\mathcal{M}}(T, \tau_2, \tau_3)}\right) \bigg/ \log\left(\frac{\tau_1}{\tau_2}\right),
\]
respectively, where $\Gamma_h^{T/\tau} = \boldsymbol{X}_h^{T/\tau}(\S^1)$, $(\tau_j, h_j)$ satisfies the relation $\tau_j = \frac{0.05}{2\pi}h_j$ for $p=2$ and $\tau_j = 0.05\left(\frac{h_j}{2\pi}\right)^{2/3}$ for $p=3$.
 \end{example}

As shown in Table~\ref{table: error of BGN-CDF}, the observed rates remain between $2.16$ and $2.22$ for the 2-stage method
and between $3.15$ and $3.16$ for the 3-stage method. These results are consistent with the designed second- and third-order convergence,
respectively.

\begin{table*}[!h]
	\centering
\begin{tabular*}{\textwidth}{@{\extracolsep{\fill}}cccccc}
		\toprule
		$\tau$ & $\mathcal{E}_{\mathcal{M}}$  & Order & $\tau$ & $\mathcal{E}_{\mathcal{M}}$  & Order \\
				\hline
				$1/800$   & $4.97\mathrm{E}{-2}$ & --- & $1/500$  & $5.65\mathrm{E}{-3}$ & --- \\
				\hline
				$1/1200$  & $2.02\mathrm{E}{-2}$ & 2.22 & $1/720$  & $1.78\mathrm{E}{-3}$ & 3.16 \\
				\hline
				$1/1600$  & $1.08\mathrm{E}{-2}$ & 2.19 & $1/980$ & $6.72\mathrm{E}{-4}$ & 3.16 \\
				\hline
				$1/2000$  & $6.65\mathrm{E}{-3}$ & 2.16 & $1/1280$ & $2.89\mathrm{E}{-4}$ & 3.15 \\
		\bottomrule
	\end{tabular*}
	\caption{\label{table: error of BGN-CDF} Errors and observed convergence orders for
\eqref{IMEX-RK-CDF: s-stage BGN+CDF} applied to CDF with the elliptical initial curve.}
\end{table*}

\begin{example}[Mesh quality]
We next examine whether the proposed high-order IMEX--RK schemes preserve the mesh-distribution properties of their underlying first-order formulations. Effective tangential motion of the mesh nodes is essential for robust simulations of geometric flows. Dziuk's method provides no tangential redistribution and may therefore lead to node accumulation during complicated evolutions. By contrast, the first-order Deckelnick--Dziuk and BGN schemes are known to produce substantially better mesh quality.

 We consider the $p=3$ versions of
\eqref{IMEX-RK-CSF: s-stage CSF+Dziuk},
\eqref{IMEX-RK-CSF: s-stage CSF+DD},
\eqref{IMEX-RK-CSF: s-stage BGN+CSF}, and
\eqref{IMEX-RK-CDF: s-stage BGN+CDF}.
As a benchmark, we use the Mikula--\v{S}ev\v{c}ovi\v{c} curve \cite{sevcovic2001evolution},
\[
\bX_0(\rho)=\left(\cos \rho,\,
\frac12\left(\sin(\cos \rho)+\sin \rho
\left(0.7+\sin \rho\,\sin^2(3\rho)\right)\right)\right),
\qquad \rho\in\mathbb{S}^1.
\]
\end{example}

\begin{figure}[h!]
	\centering
	\vspace{-3mm} \includegraphics[width=0.98\linewidth]{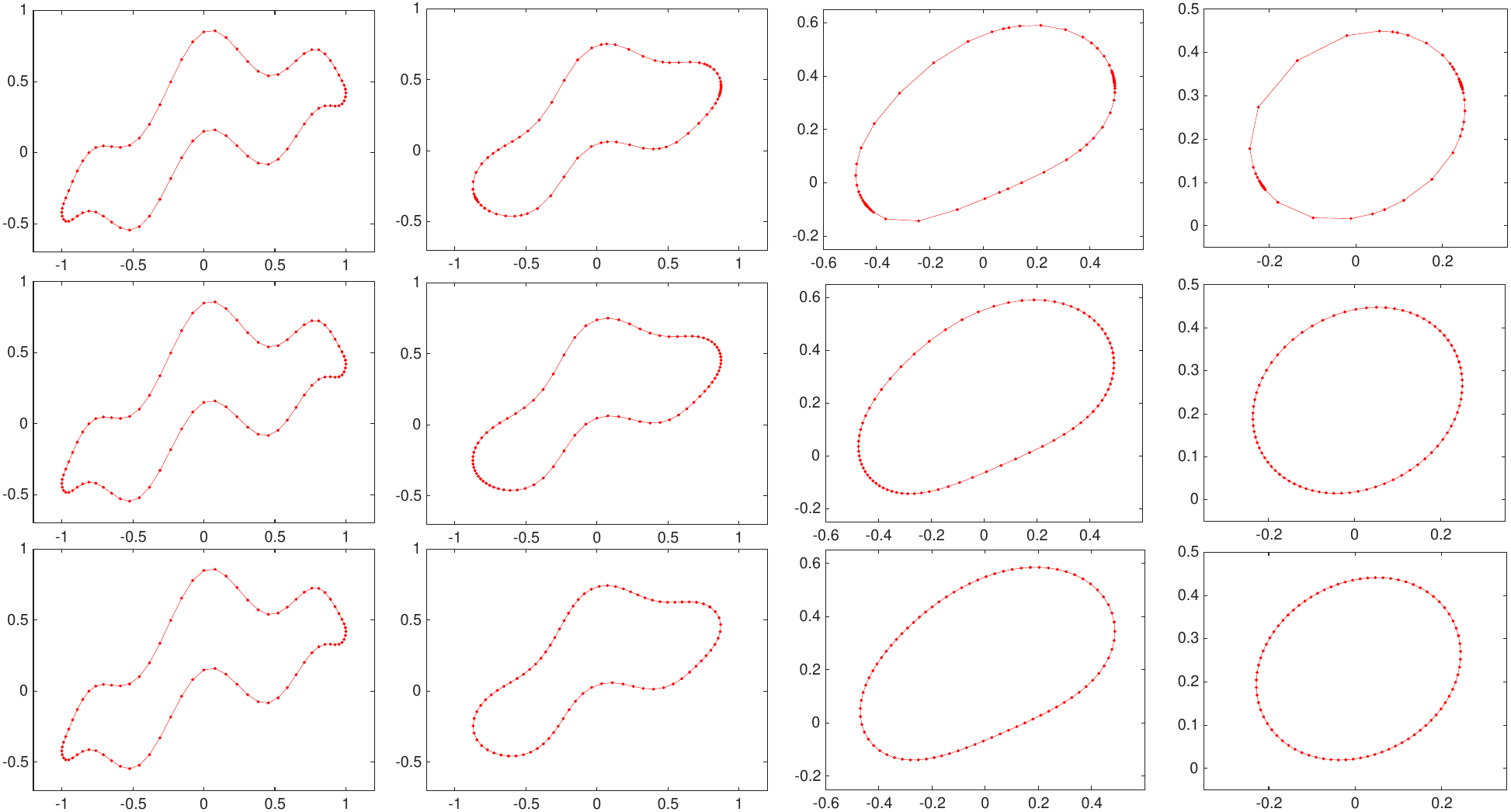}
	\vspace{-1mm}\caption{Evolution of the Mikula--\v{S}ev\v{c}ovi\v{c} curve under CSF with $p=3$, $N=80$, and $\tau=10^{-4}$. From top to bottom, the rows correspond
to the IMEX--RK discretizations based on the Dziuk, Deckelnick--Dziuk, and
BGN formulations, respectively.}
	\label{fig: Mesh_comp}
\end{figure}

Figure~\ref{fig: Mesh_comp} shows the CSF evolutions produced by the three
$p=3$ schemes. The results indicate that the IMEX--RK discretizations retain
the characteristic mesh behavior of their underlying first-order methods.
In particular, the BGN-based scheme yields the most uniform node
distribution, while the Dziuk- and Deckelnick--Dziuk-based schemes exhibit
more visible node clustering during the evolution.

\begin{figure}[h!]
	\centering	\includegraphics[width=1\linewidth]{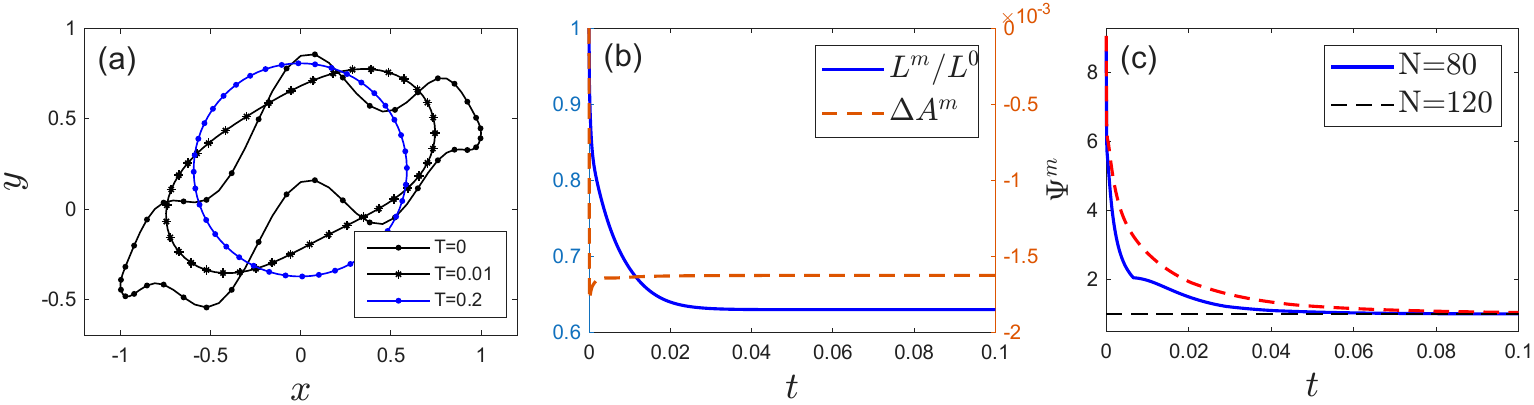}
	\vspace{-3mm}
	\caption{Mesh and geometric quantities for
\eqref{IMEX-RK-CDF: s-stage BGN+CDF} with $p=3$ and $\tau=10^{-4}$.
Panels (a)--(b) use $N=80$: (a) evolution of the curve; (b) normalized
perimeter $L^m/L^0$ and relative area change $\Delta A^m$; (c) mesh ratio $\Psi^m$. }
	\label{fig: Evo_SDF}
\end{figure}

For the CDF simulation, we further monitor several geometric and mesh
quantities. Since CDF decreases the curve length while preserving the
enclosed area, we consider the normalized perimeter $L^m/L^0$, the mesh ratio and the relative area change defined by
\[\Psi^m
=\frac{\max_j |\bX_j^m-\bX_{j-1}^m|}
{\min_j |\bX_j^m-\bX_{j-1}^m|},
\qquad
\Delta A^m=\frac{A^m-A^0}{A^0}.
\]
Here, $L^m$ denotes the perimeter of the polygonal curve $\Gamma^m$, and
$A^m$ denotes the area enclosed by $\Gamma^m$.
Figure~\ref{fig: Evo_SDF} shows that the normalized perimeter decreases
monotonically, the enclosed area is well preserved, and the mesh ratio
drops rapidly toward a nearly uniform distribution. These results indicate
that the IMEX--RK extension preserves the favorable mesh-redistribution
property of the BGN formulation for CDF.

\begin{example}[Evolution under MCF]
We next examine the performance of the proposed IMEX--RK schemes for
MCF on two benchmark dumbbell surfaces. The first surface has a relatively thick waist and is parametrized by
	\begin{align}
		\boldsymbol{X}(\theta, \varphi) = \begin{pmatrix}
			\cos \varphi \\
			\left(0.6 \cos^2 \varphi + 0.4\right) \cos \theta \sin \varphi \\
			\left(0.6 \cos^2 \varphi + 0.4\right) \sin \theta \sin \varphi
		\end{pmatrix}, \quad \theta \in [0, 2\pi),\ \varphi \in [0, \pi],\label{IMEX-RK: fat dumbbell}
	\end{align}
whereas the second has a thinner waist,
	\begin{align}
		\boldsymbol{X}(\theta, \varphi) = \begin{pmatrix}
			\cos \varphi \\
			\left(0.7 \cos^2 \varphi + 0.3\right) \cos \theta \sin \varphi \\
			\left(0.7 \cos^2 \varphi + 0.3\right) \sin \theta \sin \varphi
		\end{pmatrix}, \quad \theta \in [0, 2\pi),\ \varphi \in [0, \pi]. \label{IMEX-RK: thin dumbbell}
	\end{align}
\end{example}

To monitor the geometric evolution, we consider the normalized surface area
\[
\left.
\frac{S(t)}{S(0)}
\right|_{t=t_m}
=
\frac{S^m}{S^0},
\]
where $S^m$ denotes the area of the polyhedral surface at time $t_m$.
We also use two mesh-quality indicators introduced in
\cite{jiang2024stable}. For a triangle
$\sigma_j^m$ with vertices
$\boldsymbol{q}_{j_1}^m$,
$\boldsymbol{q}_{j_2}^m$, and
$\boldsymbol{q}_{j_3}^m$, define
\[
\left.r_h(t)\right|_{t=t_m}
:=
\frac{
\displaystyle
\max_j
\max
\left\{
\|\boldsymbol{q}_{j_1}^m-\boldsymbol{q}_{j_2}^m\|,
\|\boldsymbol{q}_{j_2}^m-\boldsymbol{q}_{j_3}^m\|,
\|\boldsymbol{q}_{j_3}^m-\boldsymbol{q}_{j_1}^m\|
\right\}
}{
\displaystyle
\min_j
\min
\left\{
\|\boldsymbol{q}_{j_1}^m-\boldsymbol{q}_{j_2}^m\|,
\|\boldsymbol{q}_{j_2}^m-\boldsymbol{q}_{j_3}^m\|,
\|\boldsymbol{q}_{j_3}^m-\boldsymbol{q}_{j_1}^m\|
\right\}
},
\]
and
\[
\left.r_a(t)\right|_{t=t_m}
:=
\frac{\max_j|\sigma_j^m|}
{\min_j|\sigma_j^m|}.
\]
Both ratios are at least one, and values closer to one indicate a more
uniform mesh.

The fat and thin dumbbells are discretized by triangular meshes with
$(J,K)=(3676,1840)$ and $(J,K)=(3276,1640)$, respectively, where $J$
and $K$ denote the numbers of triangles and vertices. We first use the fat
dumbbell \eqref{IMEX-RK: fat dumbbell} to investigate how the RK coefficients
affect tangential mesh redistribution in
\eqref{IMEX-RK-MCF: s-stage BGN+MCF}.

We begin with three $A$-stable 2-stage coefficient sets, which are shown in Table \ref{2stage2}.  Particularly, coefficient set (I) corresponds to the
predictor--corrector (PC) method.
\begin{table*}[!h]
	\centering
\def\arraystretch{1.5}
\begin{tabular*}{\textwidth}{@{\extracolsep{\fill}}ccccccccccc}
		\toprule
		& $\alpha_{10}$ & $\alpha_{11}$ & $\beta_{00}^{(1)}$ & $\beta_{10}^{(1)}$ & $\beta_{00}^{(2)}$ & $\beta_{01}^{(2)}$ & $\beta_{10}^{(2)}$  & $\beta_{11}^{(2)}$ & $\beta_{20}^{(2)}$ & $\beta_{21}^{(2)}$\\
		\hline
		(I) & $1$ & $0$ & $0$ & $\frac12$ & 0 & $\frac12$ & 0 & 0 & 0 & $\frac12 $\\
		\hline
		(II) & $\frac34$ & $\frac14$ & 0 & $2-\sqrt{2}$ & 0 & $\frac{2+\sqrt2}{4}$ &0 & $-1$ & 0 & 1\\
		\hline
		(III) & $\frac12$  & $\frac12$ & $1-\sqrt2$ & $\sqrt2$ & 0& $ \frac{\sqrt2} {2}$ & 0 & $-\frac{3+\sqrt2}{2}$& 0 & 2
		\\
		\bottomrule
	\end{tabular*}
\caption{Three sets of  coefficients for $A$-stable 2-stage methods.}
\label{2stage2}
\end{table*}

\begin{figure}[h!]
	\centering
	\includegraphics[width=1\linewidth]{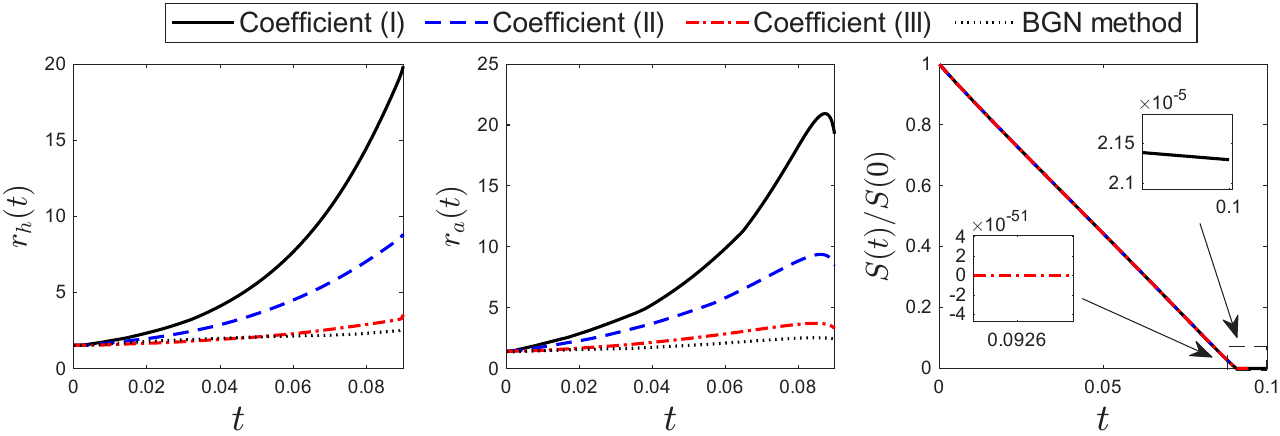}
	\vspace{-3mm}
	\caption{Evolution of geometric and mesh quantities for three 2-stage coefficient
sets and the classical BGN method. Left: edge-length ratio $r_h(t)$;
middle: triangle-area ratio $r_a(t)$; right: normalized surface area
$S(t)/S(0)$. The time step is $\tau=10^{-4}$.}
	\label{fig: RK2_3coes_ra_rh}
\end{figure}

Figure~\ref{fig: RK2_3coes_ra_rh} compares the three 2-stage coefficient sets. The mesh-quality indicators exhibit a clear ordering. For the PC coefficients (I), both $r_h(t)$ and $r_a(t)$ grow rapidly. Set (II) substantially reduces this deterioration, while set (III) gives the best mesh distribution among the three 2-stage choices and remains much closer to the classical BGN method. These observations suggest that the coefficients multiplying the stage unknowns, in particular $\beta_{10}^{(1)}$ and $\beta_{21}^{(2)}$, have a strong influence on tangential mesh redistribution.

The normalized surface-area curves are almost indistinguishable over most of the evolution, but the insets in Figure~\ref{fig: RK2_3coes_ra_rh} reveal important differences near extinction. The PC method leaves a residual normalized area of order $10^{-5}$ near $t=0.1$. By contrast, the better-performing coefficient choices drive the numerical surface area to values that are effectively zero at the displayed scale. Thus, an $A$-stable coefficient set alone is not sufficient to guarantee satisfactory mesh behavior; the particular stage coefficients also play an important role.

\begin{figure}[htbp!]
	\centering \includegraphics[width=0.9\linewidth]{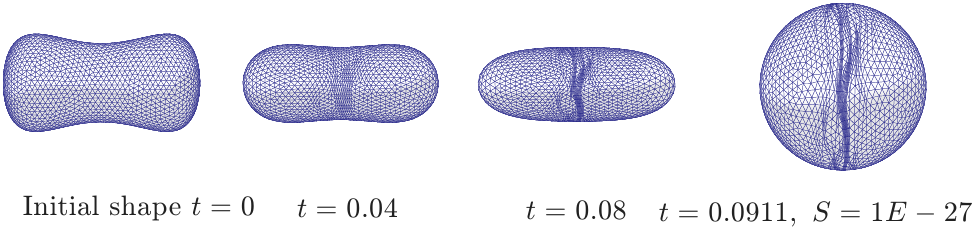}
	\vspace{-3mm}
	\caption{Evolution of the fat dumbbell \eqref{IMEX-RK: fat dumbbell} computed by
the 3-stage method \eqref{IMEX-RK-MCF: s-stage BGN+MCF} with the
coefficients in \eqref{IMEX-RK: Coefficients 3-stage}.
The time step is $\tau=10^{-4}$. The snapshots are scaled for visualization.}
	\label{fig: dumbbell6_RK3_bad}
\end{figure}

Figure~\ref{fig: dumbbell6_RK3_bad} shows the evolution obtained with the
3-stage method using the coefficients in
\eqref{IMEX-RK: Coefficients 3-stage}. The mesh remains reasonably regular
during most of the evolution. Noticeable deformation appears only in the
final rescaled snapshot, when the surface area has decreased to approximately
$10^{-27}$ and the solution is already very close to the extinction
singularity. It is therefore unclear whether this late-stage mesh distortion
is caused primarily by the RK coefficients or by the near-singular
degeneration of the surface. Nevertheless, this behavior motivates us to
seek an alternative 3-stage coefficient set with more robust mesh
redistribution near extinction.

Motivated by the preceding experiments, we employ the following $A$-stable
3-stage coefficients:
\begin{align}
\begin{split}
\alpha_{10} &= 0.77073838491,
\quad
\alpha_{11} = 0.22926161509,
\\
\alpha_{20} &= 0.99795401240,
\quad
\alpha_{21} = 0.14624530245,
\quad
\alpha_{22} = -0.14419931485,
\\
\beta_{00}^{(1)}
&= 0.40261401487,
\quad
\beta_{10}^{(1)}
= 3.02938074085,
\quad
\beta_{01}^{(2)}
= -0.10736647563,
\\
\beta_{11}^{(2)}
&= -1.54093280715,
\quad
\beta_{21}^{(2)}
= 1.72627513564,
\quad
\beta_{02}^{(3)}
= 0,
\\
\beta_{12}^{(3)}
&= -0.80613116875,
\quad
\beta_{22}^{(3)}
= -0.84039636414,
\quad
\beta_{32}^{(3)}
= 2.26931806337.
\end{split}
\label{IMEX-RK: Coefficients 3-stage set 3}
\end{align}
Compared to \eqref{IMEX-RK: Coefficients 3-stage}, this choice has relatively large values of
$\beta_{10}^{(1)}$, $\beta_{21}^{(2)}$, and $\beta_{32}^{(3)}$, motivated
by the favorable mesh behavior observed above.

\begin{figure}[htbp!]
	\centering \includegraphics[width=0.95\linewidth]{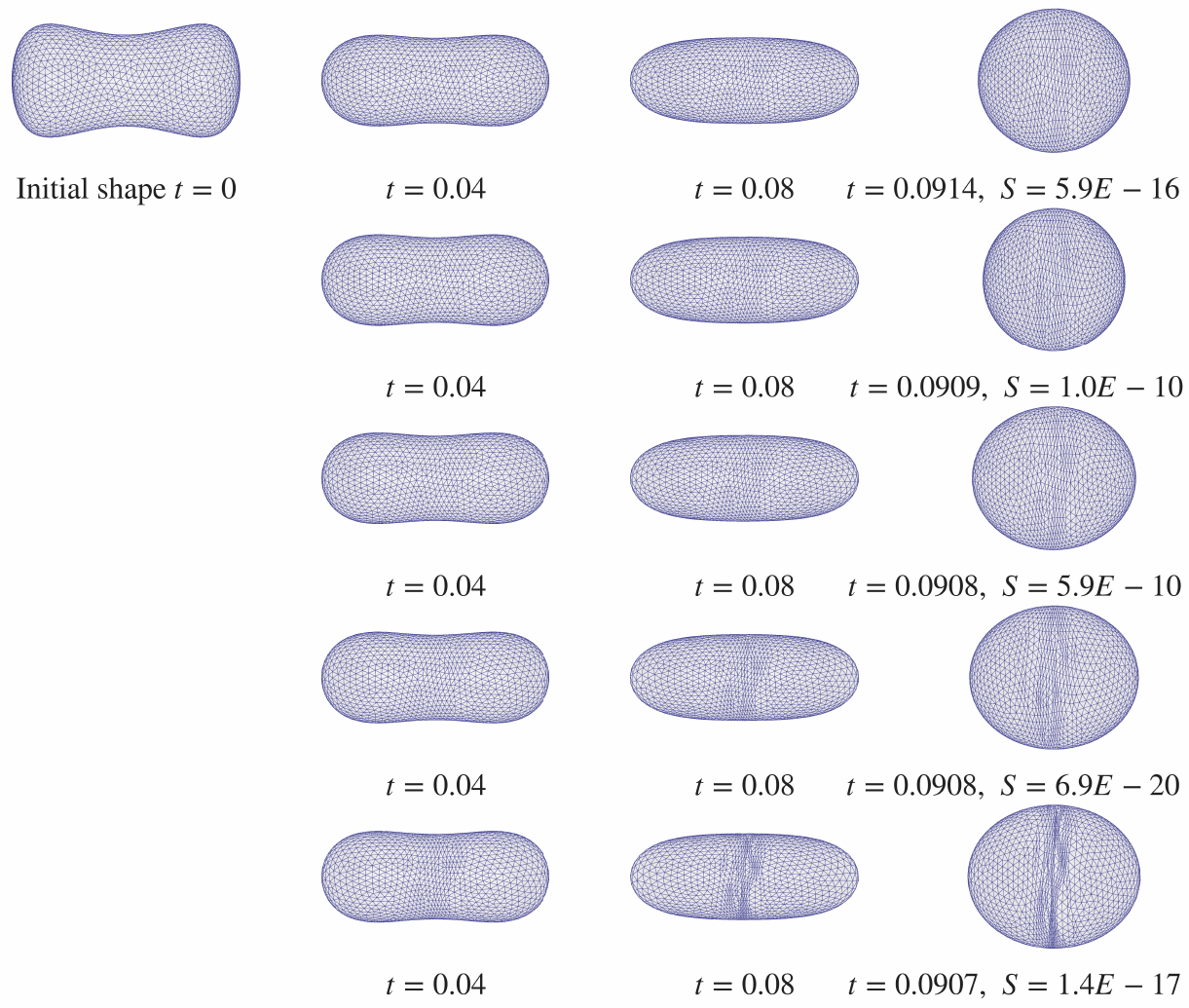}
	\vspace{-3mm}
	\caption{Evolution of the fat dumbbell \eqref{IMEX-RK: fat dumbbell}.
From top to bottom: classical BGN method;
\eqref{IMEX-RK-MCF: s-stage BGN+MCF} with $p=2$ and coefficients (III) in Table \ref{2stage2};
\eqref{IMEX-RK-MCF: s-stage BGN+MCF} with $p=3$ and coefficients
\eqref{IMEX-RK: Coefficients 3-stage set 3};
BGN/BDF2; and BGN/BDF3.
The time step is $\tau=10^{-4}$ in all simulations. The snapshots are scaled for visualization.}
	\label{fig: dumbbell6}
\end{figure}

\begin{figure}[htbp!]
	\centering
	\includegraphics[width=0.95\linewidth]{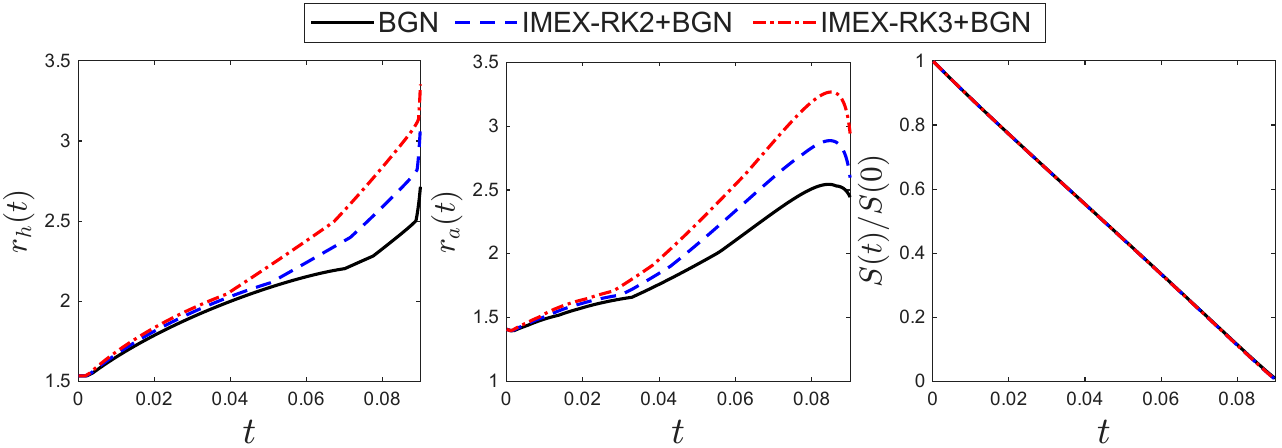}
	\vspace{-3mm}
	\caption{Evolution of geometric quantities.  Left: The mesh distribution $r_h(t)$; Middle: $r_a(t)$; Right: normalized surface area $S(t)/S(0)$. The time step is chosen as $\tau = 10^{-4}$.}
	\label{fig: rh_ra_area}
\end{figure}

\begin{figure}[htbp!]
	\centering \includegraphics[width=0.9\linewidth]{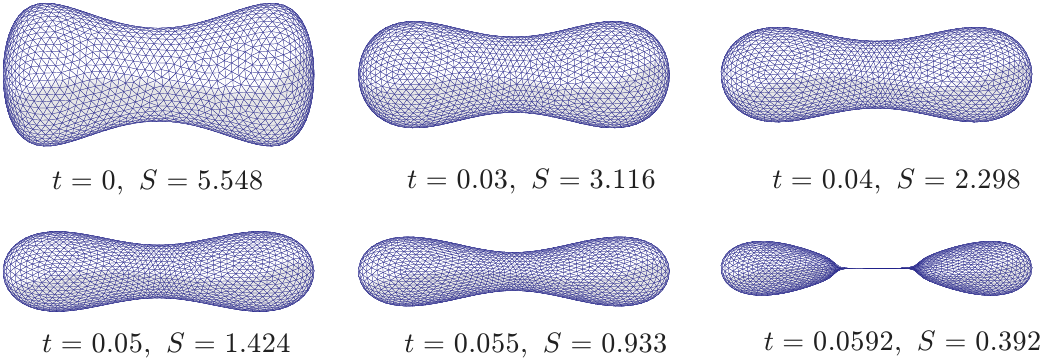}
	\vspace{-1mm}
	\caption{Evolution of the thin dumbbell \eqref{IMEX-RK: thin dumbbell} computed by
the 2-stage method \eqref{IMEX-RK-MCF: s-stage BGN+MCF} with coefficients (II) in
Table \ref{2stage2}.
The time step is $\tau=10^{-4}$. The snapshots are scaled for visualization.}
	\label{fig: Dumbbell7}
\end{figure}

Figure~\ref{fig: dumbbell6} compares five temporal discretizations for the
fat dumbbell: the classical first-order BGN method, the proposed
IMEX--RK2 and IMEX--RK3 schemes, and the BGN/BDF2 and BGN/BDF3 schemes
from \cite{jiang2024stable}. Each computation is continued until the
surface becomes nearly singular and the numerical solution can no longer
be advanced reliably.

Several observations can be made from Figure~\ref{fig: dumbbell6}.

\begin{enumerate}
\item[(1)]
The coefficient study supports the empirical observation that sufficiently
large coefficients associated with the new stage variables can enhance
tangential node redistribution. With the selected coefficients, both the
2-stage and 3-stage IMEX--RK schemes maintain regular meshes until very
close to extinction.

\item[(2)]
In this experiment, the meshes generated by the proposed IMEX--RK schemes
remain more regular near extinction than those produced by the corresponding
BGN/BDF$k$ schemes. One possible contributing factor is that the RK
construction obtains the required geometric information directly from its
intermediate stages and therefore does not require a separate prediction
procedure on a refined temporal grid.

\item[(3)]
The final resolved times for the IMEX--RK2 and IMEX--RK3 schemes are
approximately $0.0909$ and $0.0908$, respectively. These values agree
closely with the reference blow-up time
\(
T_{\mathrm{blow\text{-}up}}\approx0.0908
\) reported in \cite{m2017approximations,jiang2024stable}.
\end{enumerate}

Figure~\ref{fig: rh_ra_area} provides a quantitative comparison of the
classical BGN, IMEX--RK2, and IMEX--RK3 schemes up to $T=0.09$.
The classical BGN method gives the smallest values of $r_h(t)$ and $r_a(t)$,
while the two high-order IMEX--RK schemes produce slightly larger but still
well-controlled mesh ratios. The normalized surface-area curves are nearly
indistinguishable. Hence, the proposed high-order schemes retain satisfactory
mesh quality while accurately reproducing the rapid area decay as the
surface approaches extinction.

For the thin dumbbell \eqref{IMEX-RK: thin dumbbell}, we apply the 2-stage method \eqref{IMEX-RK-MCF: s-stage BGN+MCF} with coefficients (II) in Table \ref{2stage2}. Figure~\ref{fig: Dumbbell7} shows the formation of a neck-pinch singularity at approximately $t=0.0592$, in agreement with the reference value reported in \cite{jiang2024stable}.

\begin{example}[Evolution under SDF]
We finally consider SDF for two initial surfaces: an ellipsoid
\[
x^2+4y^2+4z^2=4,
\]
and an $8\times1\times1$ cuboid.
\end{example}

\begin{figure}[htbp!]
	\centering \includegraphics[width=0.95\linewidth]{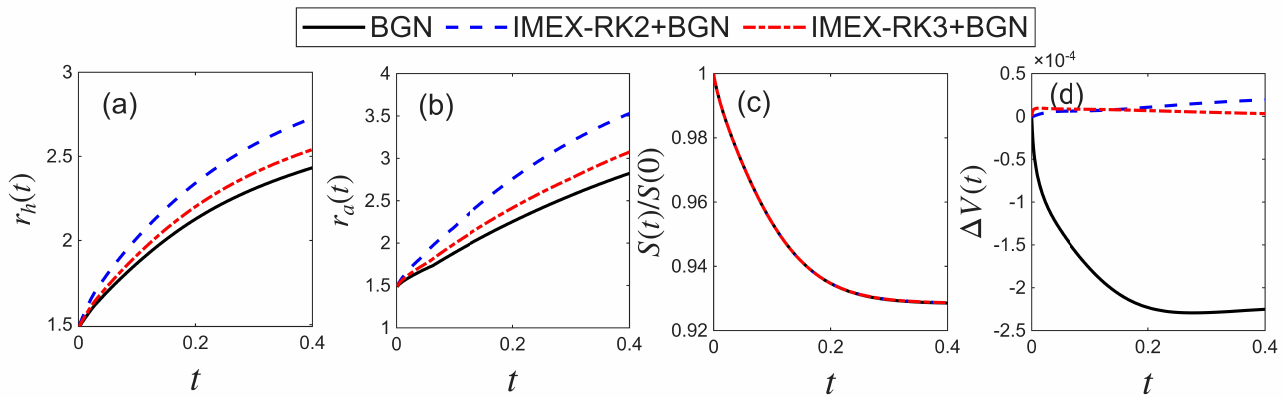}
	\vspace{-3mm}
	\caption{Evolution of geometric quantities for the ellipsoid under SDF: (a) edge-length ratio $r_h(t)$; (b) triangle-area ratio $r_a(t)$; (c) normalized surface area $S(t)/S(0)$; and (d) relative volume change $\Delta V(t)$. The time step is $\tau=10^{-4}$.}
	\label{fig: rh_ra_area_volume}
\end{figure}

\begin{figure}[htbp!]
	\centering
	\includegraphics[width=0.95\linewidth]{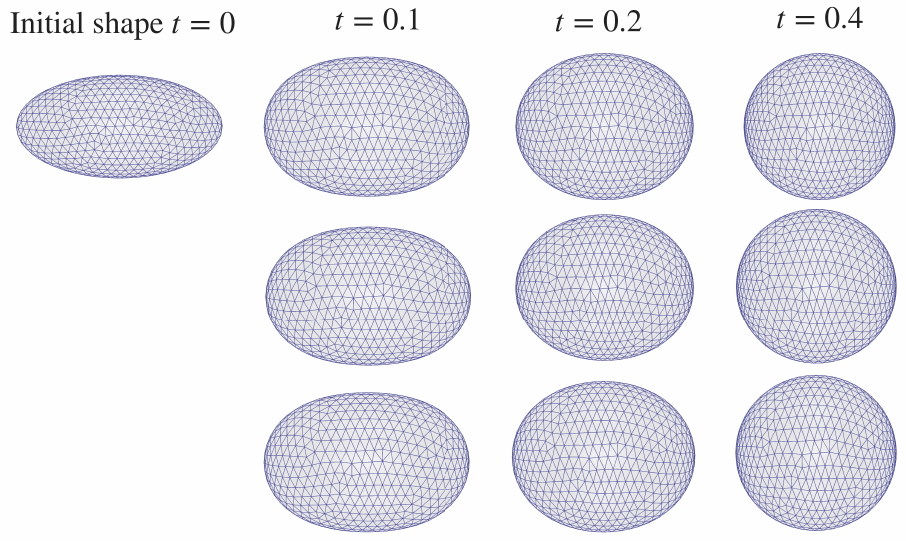}
	\vspace{-3mm}
	\caption{Evolution of the ellipsoid under SDF. First row: classical BGN method. Second row: \eqref{IMEX-RK-SDF: s-stage BGN+SDF} with $p=2$ and coefficients (III) in Table \ref{2stage2}. Third row: \eqref{IMEX-RK-SDF: s-stage BGN+SDF} with $p=3$ and coefficients \eqref{IMEX-RK: Coefficients 3-stage set 3}.}
	\label{fig: SDF_ellipsoid}
\end{figure}

\begin{figure}[htbp!]
	\centering \includegraphics[width=0.95\linewidth]{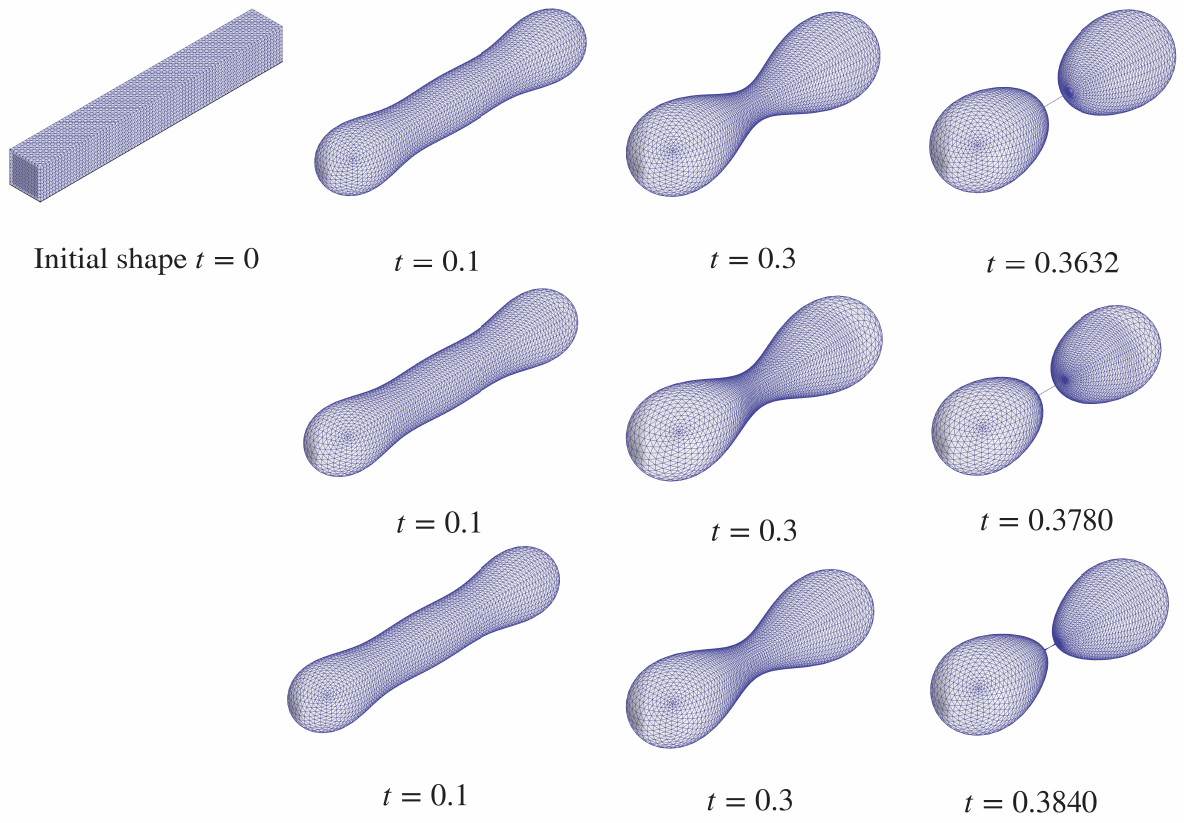}
	\vspace{-3mm}
	\caption{Evolution of the cuboid under SDF. First row: classical BGN method. Second row: \eqref{IMEX-RK-SDF: s-stage BGN+SDF} with $p=2$ and coefficients (III) in Table \ref{2stage2}. Third row: \eqref{IMEX-RK-SDF: s-stage BGN+SDF} with $p=3$ and coefficients \eqref{IMEX-RK: Coefficients 3-stage set 3}. }
	\label{fig: SDF_cuboid}
\end{figure}

For the ellipsoid, we monitor $r_h(t)$, $r_a(t)$, the normalized surface area $S(t)/S(0)$, and the relative change in enclosed volume
\[
\Delta V(t)
:=
\frac{V(t)-V(0)}{V(0)}.
\]
The results are shown in Figure~\ref{fig: rh_ra_area_volume}. All three methods reproduce the expected geometric behavior of SDF: the surface area decreases monotonically, while the enclosed volume is approximately conserved. The magnitude of the relative volume error produced by the IMEX--RK schemes is substantially smaller than that of the classical BGN method in this test. Figure~\ref{fig: SDF_ellipsoid} shows the corresponding surface evolutions for $(J,K)=(1556,780)$. The meshes remain regular throughout the displayed evolution, demonstrating the robustness of the proposed high-order schemes.

For the cuboid, we use $(J,K)=(6800,3402)$. Figure~\ref{fig: SDF_cuboid} shows the evolution toward pinch-off. At the spatial and temporal resolutions used here, the detected pinch-off times are approximately
\[
t=0.3632
\quad\text{for BGN};\quad
t=0.3780
\quad\text{for IMEX--RK2};\quad
t=0.3840
\quad\text{for IMEX--RK3}.
\]
For this benchmark, Barrett et al. reported $t\approx0.369$ using an adaptive-in-time computation \cite[Section~5.5]{barrett2008parametric}. A value of approximately $0.370$ was subsequently reported in \cite[Section~4.2]{bao2021structure} using an implicit method and in \cite[Section~4.2]{jiang2024stable} using BGN/BDF$k$ schemes. The pinch-off times obtained with the present IMEX--RK discretizations are therefore somewhat later than the previously reported values. A systematic space--time refinement study would be needed to determine the asymptotic accuracy of the detected pinch-off time.

\section{Conclusion}

In this paper, we developed a high-order IMEX--RK framework for parametric finite element approximations of geometric flows. Starting from Dziuk's semidiscrete formulation of curve shortening flow, we recast the semidiscrete equations as a finite-dimensional ODE system and constructed second- and third-order IMEX--RK schemes. Each stage requires only a linear solve, and the intermediate stages generate the geometric predictions needed at subsequent stages, eliminating the need for a separate prediction procedure. The one-step structure also avoids the additional starting values required by multistep methods. We further identified $A$-stability as an important criterion for guiding the selection of RK coefficients. The same construction extends naturally to BGN formulations, with the geometry generated at the preceding stage serving as both the integration domain and the predictor geometry for the next stage. This extension provides a unified framework for high-order approximations of mean curvature flow and surface diffusion for curves and surfaces. Numerical experiments support the designed orders of temporal convergence and indicate that, with suitable coefficient choices, the proposed schemes maintain mesh quality comparable to that of the underlying first-order methods. The resulting framework combines high-order temporal accuracy, linear stage solves, internally generated geometric predictions, and effective mesh redistribution.


\end{document}